\documentclass[A4,11pt]{article}
\usepackage{amsfonts}
\usepackage{amsmath,amsthm}
\usepackage{amssymb,amscd}
\usepackage{blkarray}
\usepackage{enumerate}
\usepackage{kbordermatrix}
\usepackage{color}
\usepackage{cite}
\usepackage{esint}
\usepackage{tikz}
\usepackage{tikzscale}
\usepackage[absolute,overlay]{textpos}
\usepackage{multicol}
\usepackage{framed}
\usepackage{lscape}
\usepackage{caption}
\usepackage{subcaption}
\usepackage[inline]{enumitem}
\usepackage{float}
\usepackage{mathtools}
\usepackage{lipsum}
\renewcommand{\kbldelim}{(}% Left delimiter
\renewcommand{\kbrdelim}{)}% Right delimiter
\makeatletter

\newcommand{\Rmnum}[1]{\expandafter\@slowromancap\romannumeral #1@}
\makeatother

\def\XXint#1#2#3{{\setbox0=\hbox{$#1{#2#3}{\int}$}
\vcenter{\hbox{$#2#3$}}\kern-.5\wd0}}

\def\XXint#1#2#3{{\setbox0=\hbox{$#1{#2#3}{\int}$}
\vcenter{\hbox{$#2#3$}}\kern-.5\wd0}}

\theoremstyle{definition}
\newtheorem{theorem}{Theorem}[section]

\newtheorem{lemma}[theorem]{Lemma}

\newtheorem{definition}{Definition}[section]
\newtheorem{proposition}[theorem]{Proposition}

\newtheorem{remark}{Remark}[section]

\begin{document}
\title{Integrability of the multi-species ASEP with long-range jumps on $\mathbb{Z}$}
\author{\textbf{Eunghyun Lee\footnote{eunghyun.lee@nu.edu.kz}}\\ {\text{Department of Mathematics,}}
                                         \date{}   \\ {\text{School of Sciences and Humanities,}} \\ {\text{Nazarbayev University, }}\\ {\text{Kazakhstan }}   }

\date{}
\maketitle
\begin{abstract}
\noindent Let us consider a two-sided multi-species stochastic particle model with finitely many particles on $\mathbb{Z}$ defined as follows. Suppose that each particle is labelled by a positive integer $l$ and waits a random time exponentially distributed with rate $1$. It then chooses the right direction to jump with probability $p$ or the left direction with probability $q=1-p$. If the particle chooses the right direction, it jumps to the nearest site occupied by a particle $l'<l$ (with the convention that an empty site is considered as a particle with labelled $0$). If the particle chooses the left direction, it follows the rule of the multi-species totally asymmetric simple exclusion process (mTASEP). We show that this model is integrable and provide the exact formula of the transition probability using the Bethe ansatz. 
\end{abstract}

\section{Introduction}
In this paper, we treat a multi-species stochastic particle model with finitely many particles on the infinite lattice $\mathbb{Z}$, which has the exclusion rule that each site cannot be occupied by more than one particle. Each particle is labelled by a positive integer $l$ representing its species. A state of the system  can be described by specifying each particle's species and its location.  To be more specific, if  $X=(x_1,\dots, x_N)\in \mathbb{Z}^N$ with $x_1<\dots<x_N$ represents the positions of $N$ particles and the species of the particle at $x_i$ is denoted by $\pi_i$, then we denote the state of the system by $(X,\pi)=(x_1,\dots, x_N, \pi_1\cdots\pi_N)$ with $\pi = \pi_1\cdots \pi_N$. In a single-species model, of course, a state is simply denoted by $(x_1,\dots, x_N)$.

The asymmetric simple exclusion process (ASEP) is the most well-known single-species stochastic particle model. In the ASEP on $\mathbb{Z}$, each particle waits an exponential time with rate 1 to choose a direction to jump. It chooses the right direction with probability $p$ or the left direction with probability $q=1-p$. After choosing a direction, a particle at $x$ immediately jumps to $x\pm 1$ if the  site is empty. However, if the target site is already occupied by another particle, the particle gives up jumping and waits a random time again to jump. The technique to find the transition probability of the ASEP on $\mathbb{Z}$ using the Bethe ansatz is well known \cite{Schutz-1997,Tracy-Widom-2008}.

In the multi-species ASEP, a higher-numbered particle has priority over a lower-numbered particle in jumping. To be more specific, if a particle $l$ at $x$ chooses to jump to $x \pm 1$  already occupied by $l'<l$, then the particle $l$ can indeed jump to the site by switching their positions.  However, if $l'\geq l$, then the particle  $l$ gives up jumping to the site occupied by $l'$. This model has been extensively studied in many directions \cite{Chatterjee-Hayakawa,Chatterjee-Schutz-2010,Chen-Gier-Hiki-Sasamoto-Usui,Gier-Mead-Wheeler-2021,Kuniba-Maruyama-Okado,Kuan-2020,Lee-2017,Lee-2018,Prolhac-Evans-Mallick-2009,Tracy-Widom-2009,Tracy-Widom-2013}.

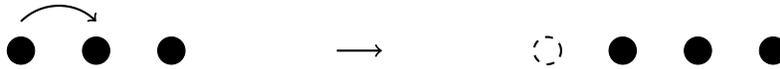
\begin{figure}[h]
  \centering
  \begin{tikzpicture}
    % circles with numbers inside
    \node[circle, draw, fill=black, minimum size=0.3cm] (left1) at (0,0) {};
    \node[circle, draw, fill=black, minimum size=0.3cm] (left2) at (1,0) {};
    \node[circle, draw, fill=black, minimum size=0.3cm] (left3) at (2,0) {};
    \node[circle, minimum size=0.3cm] (left4) at (3,0) {};

    % curved arrow
    \draw[->, thick, bend left=45] ([yshift=0.2cm]left1.north) to ([yshift=0.2cm]left2.north);

    % Arrow pointing to the right
    \draw[->, thick] (4.2,0) -- (4.8,0);

    % circles with numbers inside
    \node[circle, draw, dashed, thick, minimum size=0.3cm] at (7,0) {};
    \node[circle, draw, fill=black, minimum size=0.3cm] at (8,0) {};
    \node[circle, draw, fill=black, minimum size=0.3cm] at (9,0) {};
    \node[circle, draw, fill=black, minimum size=0.3cm] at (10,0) {};

  \end{tikzpicture}
   \caption{Drop-push model}
  \label{drop-push}
\end{figure}

Another stochastic particle model with rules as simple as those of the ASEP is the drop-push model \cite{Ali,Sasamoto-Wadati-1998,Schutz-Ramaswamy-Barma}. In the drop-push model, each particle jumps to the next site on the right after an exponential random time with rate $1$. Unlike the ASEP, a particle can jump even if the target site is already occupied by pushing the particle occupying the site. For example, if  a particle at $x$ jumps to $x+1$ occupied by another particle, the particle at $x$ jumps to $x+1$, pushing the particle at $x+1$ to $x+2$. In this case, if $x+2$ is also occupied, then the particle pushed to $x+2$ also pushes the particle at $x+2$ to $x+3$ (Figure \ref{drop-push} illustrates the definition of the drop-push model).

\begin{figure}[h]
    \centering
    \begin{tikzpicture}
        % Left particles
        \node[circle, draw, fill=black, minimum size=0.5cm] (left1) at (0,2) {};
        \node[circle, draw, fill=white, minimum size=0.5cm] (left2) at (1,2) {};
        \node[circle, draw, fill=black, minimum size=0.5cm] (left3) at (2,2) {};
        \node[circle, draw, fill=white, minimum size=0.5cm] (left4) at (3,2) {};
        \node[circle, draw, fill=black, minimum size=0.5cm] (left5) at (4,2) {};

        \node[circle, draw, fill=black, minimum size=0.5cm] at (0,1) {};
        \node[circle, draw, fill=black, minimum size=0.5cm] at (1,1) {};
        \node[circle, draw, fill=black, minimum size=0.5cm] at (2,1) {};
        \node[circle, draw, fill=black, minimum size=0.5cm] at (3,1) {};
        \node[circle, draw, fill=black, minimum size=0.5cm] at (4,1) {};

        % Arrow pointing to the right
        \draw[->,thick] (5.8,0) -- (6.3,0);

        % Curved arrow from left3 to left2
        \draw[<-, thick, bend left=45] ([yshift=0.2cm]left2.north) to ([yshift=0.2cm]left3.north);

        % Curved arrow from left4 to left5
        \draw[->, thick, bend left=45] ([yshift=0.2cm]left4.north) to ([yshift=0.2cm]left5.north);

        % Right particles
        \node[circle, draw, fill=black, minimum size=0.5cm] at (7,2) {};
        \node[circle, draw, fill=black, minimum size=0.5cm] at (8,2) {};
        \node[circle, draw, fill=white, minimum size=0.5cm] at (9,2) {};
        \node[circle, draw, fill=white, minimum size=0.5cm] at (10,2) {};
        \node[circle, draw, fill=black, minimum size=0.5cm] at (11,2) {};
        \node[circle, draw, fill=white, minimum size=0.5cm] at (12,2) {};

        \node[circle, draw, fill=black, minimum size=0.5cm] at (7,1) {};
        \node[circle, draw, fill=black, minimum size=0.5cm] at (8,1) {};
        \node[circle, draw, fill=black, minimum size=0.5cm] at (9,1) {};
        \node[circle, draw, fill=white, minimum size=0.5cm] at (10,1) {};
        \node[circle, draw, fill=black, minimum size=0.5cm] at (11,1) {};
        \node[circle, draw, fill=black, minimum size=0.5cm] at (12,1) {};
    \end{tikzpicture}

    \begin{tikzpicture}
        % Left particles
        \node[circle, draw, minimum size=0.5cm, inner sep=0pt] (left1) at (0,0) {2};
        \node[circle, draw, minimum size=0.5cm, inner sep=0pt] (left2) at (1,0) {1};
        \node[circle, draw, minimum size=0.5cm, inner sep=0pt] (left3) at (2,0) {2};
        \node[circle, draw, minimum size=0.5cm, inner sep=0pt] (left4) at (3,0) {1};
        \node[circle, draw, minimum size=0.5cm, inner sep=0pt] (left5) at (4,0) {2};
        \node[circle, draw=none, minimum size=0.5cm, inner sep=0pt] (left6) at (5,0) {};

        % Curved arrow from left3 to left2
        \draw[<-, thick, bend left=45] ([yshift=0.2cm]left2.north) to ([yshift=0.2cm]left3.north);

        % Curved arrow from left2 to left3
        \draw[->, thick, bend right=45] ([yshift=-0.2cm]left2.south) to ([yshift=-0.2cm]left3.south);

        % Curved arrow from left4 to left6(empty space)
        \draw[->, thick, bend left=45] ([yshift=0.2cm]left4.north) to ([yshift=0.2cm]left6.north);

        % Right particles
        \node[circle, draw, minimum size=0.5cm, inner sep=0pt] (right1) at (7,0) {2};
        \node[circle, draw, minimum size=0.5cm, inner sep=0pt] (right2) at (8,0) {2};
        \node[circle, draw, minimum size=0.5cm, inner sep=0pt] (right3) at (9,0) {1};
        \node[circle, draw, minimum size=0.5cm, inner sep=0pt] (right4) at (10,0) {};
        \node[circle, draw, minimum size=0.5cm, inner sep=0pt] (right5) at (11,0) {2};
        \node[circle, draw, minimum size=0.5cm, inner sep=0pt] (right6) at (12,0) {1};

    \end{tikzpicture}
    \caption{Basic coupling}
    \label{fig:one25}
\end{figure}
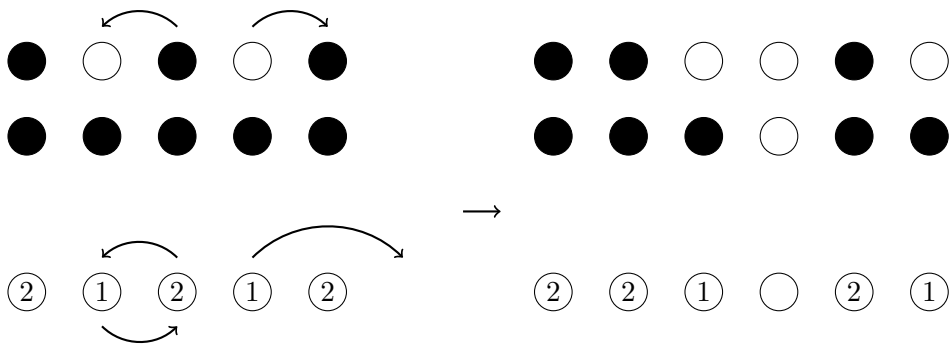

 A multi-species version of the (one-sided) drop-push model was introduced as the \textit{frog model} in \cite{Bukh-Cox}. Additionally, the stationary distribution of the \textit{frog model}  was recently studied \cite{Aggarwal-Nicoletti-Petrov,Ayyer-Martin}. In particular, \cite{Aggarwal-Nicoletti-Petrov} treated various multi-species particle models in a unified way to study the stationary distribution. Another direction of extending to a multi-species version of the (one-sided) drop-push model is found in \cite{Roshani-Khorrami}. A straightforward extension of the  one-sided drop-push model to the one where a particle jumps to the right with rate $p$ or to the left with rate $q=1-p$ is not integrable in the sense that the Bethe ansatz is not applicable. However, it is known that if the jumping rates satisfy certain conditions, the two-sided model is integrable \cite{Ali2}.  A natural question arises as to whether there exists an integrable multi-species version of the model in \cite{Ali2}. Unfortunately, the authors' best efforts to find such a version have not been fruitful. Nevertheless, in this paper, we introduce a two-sided integrable multi-species model incorporating  the rule of the long-range jumps partially.
\\ \\
\textbf{Definition of the model}\\ \\
In this paper, we show that if particles follow the rule of the frog model when they jump to the right and follow the rule of the multi-species TASEP when they jump to the left, then the corresponding two-sided model is integrable. Hence, we can use the Bethe ansatz to find the transition probability. For our task, overall, we will take a similar approach to the method in \cite{Lee-2020}. The precise definition of the model is as follows.
\begin{definition}\label{836pm72}
A particle  $l$ at $x$ waits an exponential time with rate 1 to jump to the right with probability $p$ or to the left with probability $q=1-p$. If $l$ chooses the right direction, it jumps to the nearest right site occupied by  $l'<l$ (including  $l'=0$ which corresponds to an empty site), and  $l'$ immediately jumps to the nearest right site occupied by $l''<l'$. If  $l$ at $x$ chooses the left direction to jump and $x-1$ is occupied by $l'<l$, then $l$ jumps to $x-1$ and $l'$  moves to $x$. However, if  $l'\geq l$, then $l$ gives up jumping to $x-1$.
\end{definition}
We will call this model the multi-species ASEP  with long-range jumps. The definition of our model with $p=1$ can be motivated by the basic coupling  of two drop-push models with different configurations. The idea of viewing a low-numbered particle in the two-species ASEP as a discrepancy between the configurations of two ASEPs is well known \cite{Liggett,Tracy-Widom-2009}. Using the same idea, let $\eta_t$ be a drop-push model and $\eta_t(x)$ represent the state of site $x$ at time $t$; that is, $\eta_t(x) = 0$ if $x$ is unoccupied and $\eta_t(x) = 1$ if it is occupied. Let $(\eta_t, \zeta_t)$ be the coupled process of the drop-push models $\eta_t$ and $\zeta_t$ with $\eta_0(x) \leq \zeta_0(x)$ for all $x$. If we think of a site $x$ with $\eta_t(x) = \zeta_t(x) = 1$ as being occupied by a particle of species 2, and  a site $x$ with $\eta_t(x) =0, \zeta_t(x) = 1$ as being occupied by a particle of species 1, then the evolution of the coupled process is  the same as that of our model $p=1$ in Definition \ref{836pm72}. Similarly, the two-sided model with general $p$ in Definition \ref{836pm72} can be interpreted via the basic coupling of  two-sided corresponding single species models. The two-sided single-species model where a particle follows the rule of the drop-push model when it moves to the right and follows the rule of the TASEP when it moves to the left was studied in \cite{Borodin2}, so our model can be seen as a multi-species version of the model in \cite{Borodin2}.  Figure \ref{fig:one25} illustrates the interpretation by the basic coupling.

\textbf{Notations}
\\ \\
Suppose that there are $N$ particles labelled by positive integers $a_1,\dots, a_N$, and let $\mathcal{M}=[a_1,\dots,a_N]$ be the multi-set of those $N$ integers. Let $P_{(Y,\nu)}(X,\pi;t)$ be the probability that the system will be in state $(X,\pi)$ at time $t$ given that the system is in state $(Y,\nu)$ at time $0$. We view $\pi$ and $\nu$  as permutations of  $\mathcal{M}$ of $N$ elements. Of course, if $\pi$ and $\nu$ are permutations of two different multi-sets, then  $P_{(Y,\nu)}(X,\pi;t)=0$.

We will consider an arbitrary multi-set $[a_1,\dots,a_N]$ with each $a_i \in \{1,\dots, N\}$ to cover  all possible compositions of species. For this purpose, let $\mathbf{P}_Y(X;t)$ be  the $N^N \times N^N$ matrix such that  rows and columns are labelled with $1\cdots 1, \dots, N\cdots N$ and the $(\pi,\nu)$-element is $P_{(Y,\nu)}(X,\pi;t)$.  According to the definition of the model, $\mathbf{P}_Y(X;t)$ is a lower triangular matrix when rows and columns are labelled lexicographically from $1 \cdots 1$ to $N \cdots N$.  For example, for $N=2$,
\begin{equation}\label{426pm5224}
\mathbf{P}_{Y}(X;t) =
\kbordermatrix{
     & 11 &12 &21 &22 \\[4pt]
    11 & P_{(Y,11)}(X,11;t) &  0 &  0&  0 \\[4pt]
    12  &0 &   P_{(Y,12)}(X,12;t) & 0&   0 \\[4pt]
    21 & 0 &   P_{(Y,12)}(X,21;t) &  P_{(Y,21)}(X,21;t)&   0 \\[4pt]
    22 &0 &  0& 0&   P_{(Y,22)}(X,22;t)
    }.
\end{equation}
Throughout this paper, $\mathbf{I}$ denotes the identity matrix and $\mathbf{0}$ denotes the zero matrix. If needed, we write   $\mathbf{I}_n$ and $\mathbf{0}_n$ to specify that their dimensions are $n$ by $n$. We write  $\mathbf{M}_{(\pi,\nu)}$ for  the $(\pi,\nu)$-element of a matrix $\mathbf{M}$.
\\ \\
\textbf{Organization of the paper}\\ \\
In Section \ref{7522pm1024}, we investigate the case of $N=2$ as a preliminary building block for the case of general $N$.  In Section \ref{924pm72}, we investigate the case of $N=3$. The arguments for $N=3$  build upon the results for $N=2$. However, in the case of $N=3$, \textit{two-body reducibility} and the Yang-Baxter equation are discussed for integrability. Moving from $N=2$ to $N=3$ is not a trivial extension, but if the integrability of $N=3$ is proved, then the integrability of general $N$ is trivially confirmed. Hence, we provide detailed  arguments for the case of $N=3$ in Section  \ref{924pm72} and only the results for general $N$  in Section \ref{1055am729}. This approach also helps us understand the integrability of the model better, reducing the heavy notations for general $N$. In Section \ref{1055am729}, we provide the exact formula for the transition probability for an arbitrary $N$-particle system.
\section{$N=2$}\label{7522pm1024}
\subsection{Master equations}
Each element $P_{(Y,\nu)}(X,\pi;t)$ in (\ref{426pm5224}) satisfies the master equation  and the initial condition:
\begin{equation}\label{706am725}
P_{(Y,\nu)}(X,\pi;0) =
\begin{cases}
1 &~~\textrm{if $(Y,\nu)=(X,\pi)$}\\[5pt]
0&~~\textrm{otherwise}.
\end{cases}
\end{equation}
For a fixed $Y$, the form of the master equation for $P_{(Y,\nu)}(X,\pi;t)$ depends on $X, \pi,\nu$. The master equation for $P_{(Y,\nu)}(X,\pi;t)$ when $X=(x_1,x_2)$ with  $x_1< x_2-1$ is the $(\pi,\nu)$-element of the matrix equation:
\begin{equation}\label{748pm516}
\begin{aligned}
\frac{d}{dt}\mathbf{P}_{Y}(x_1,x_2;t) = &~ p\,\mathbf{P}_{Y}(x_1-1,x_2;t) + p\,\mathbf{P}_{Y}(x_1,x_2-1;t) \\
                                      + &~ q\,\mathbf{P}_{Y}(x_1+1,x_2;t) + q\,\mathbf{P}_{Y}(x_1,x_2+1;t) - 2\mathbf{P}_{Y}(x_1,x_2;t)
 \end{aligned}
\end{equation}
where the differentiation on the left side in  (\ref{748pm516}) denotes the differentiation of each matrix element.

The master equation for $P_{(Y,\nu)}(X,\pi;t)$  when  $X=(x,x+1)$ is
the $(\pi,\nu)$-element of the matrix equation:
\begin{equation}\label{752pm516}
\begin{aligned}
&\frac{d}{dt}\mathbf{P}_Y(x,x+1;t) = p\mathbf{P}_Y(x-1,x+1;t) +   p\mathbf{B}\,\mathbf{P}_Y(x-1,x;t) - 2p\mathbf{P}_Y(x,x+1;t) \\
 &\hspace{1cm}+ q\mathbf{P}_Y(x,x+2;t)  + q\mathbf{B}_1\,\mathbf{P}_Y(x,x+1;t)- q\mathbf{B}_2\mathbf{P}_Y(x,x+1;t)
\end{aligned}
\end{equation}
where
\begin{equation}\label{Bmatrix}
\mathbf{B} = \kbordermatrix{
    & 11 & 12 & 21 & 22  \\
    11 & 1 & 0 & 0 & 0  \\
    12 & 0 & 0 & 0 & 0  \\
    21 & 0 & 1 & 1 & 0  \\
    22 & 0 & 0 & 0 & 1
  },
\end{equation}
\begin{equation}\label{Bmatrix2}
\mathbf{B}_1 = \kbordermatrix{
    & 11 & 12 & 21 & 22  \\
    11 & 0 & 0 & 0 & 0  \\
    12 & 0 & 0 & 0 & 0  \\
    21 & 0 & 1 & 0 & 0  \\
    22 & 0 & 0 & 0 & 0
  },
\end{equation}
and
\begin{equation}\label{Bmatrix3}
 \mathbf{B}_2 = \kbordermatrix{
    & 11 & 12 & 21 & 22  \\
    11 & 1 & 0 & 0 & 0  \\
    12 & 0 & 2 & 0 & 0  \\
    21 & 0 & 0 & 1 & 0  \\
    22 & 0 & 0 & 0 & 1
  }.
\end{equation}
Note that the $(ij,lm)$-element of $p\mathbf{B}$ represents the instantaneous rate of transition from state $(x-1,x,lm)$ to state $(x,x+1, ij)$, while the $(ij,lm)$-element of $q\mathbf{B}_1$ represents the instantaneous rate of transition from state $(x,x+1,lm)$ to state $(x,x+1, ij)$. Moreover,  (\ref{752pm516}) can be expressed as:
\begin{equation}\label{752pm752}
\begin{aligned}
\frac{d}{dt}\mathbf{P}_Y(x,x+1;t) =~& p\mathbf{P}_Y(x-1,x+1;t) +   p\mathbf{B}\,\mathbf{P}_Y(x-1,x;t)  \\
 &~\,+ q\mathbf{P}_Y(x,x+2;t)  + q\mathbf{B}\,\mathbf{P}_Y(x,x+1;t)- 2\mathbf{P}_Y(x,x+1;t)
\end{aligned}
\end{equation}
due to the relation $\mathbf{B} - \mathbf{B}_1 +  \mathbf{B}_2  = 2\mathbf{I}$.

The matrix form of the initial condition (\ref{706am725})  is given by:
\begin{equation}\label{654pm72}
\mathbf{P}_Y(X;0) =
\begin{cases}
\mathbf{I} &~\textrm{if ${X} = Y$}, \\
\mathbf{0} &~\textrm{if ${X} \neq {Y}$}.
\end{cases}
\end{equation}
\subsection{Bethe ansatz}\label{646pm624}
In this section, we extend the  Bethe ansatz method \cite{Lee-2020,Schutz-1997,Tracy-Widom-2008} to find solutions of  the matrix equations  (\ref{748pm516}) and (\ref{752pm752}).

Recall  that $\mathbf{P}_Y(X;t)$ in (\ref{426pm5224}) was defined for the \textit{physical region}, meaning for $Y=(y_1,y_2)$ and $X=(x_1,x_2)$ with $y_1<y_2$ and $x_1<x_2$, respectively. Let $\mathbf{U}(X;t)$ be the $4\times 4$ matrix obtained by replacing $P_{(Y,\nu)}(X,\pi;t)$ in (\ref{426pm5224}) with some functions $U_{(Y,\nu)}(X,\pi;t)$ defined for all  $X,Y \in \mathbb{Z}^2$. Suppose that $\mathbf{U}(X;t)= \mathbf{U}(x_1,x_2;t)$ satisfies
\begin{equation}\label{748pm51643}
\begin{aligned}
\frac{d}{dt}\mathbf{U}(x_1,x_2;t) = ~& p\mathbf{U}(x_1-1,x_2;t) + p\mathbf{U}(x_1,x_2-1;t) \\
+~& q\mathbf{U}(x_1+1,x_2;t) + q\mathbf{U}(x_1,x_2+1;t) - 2\mathbf{U}(x_1,x_2;t)
\end{aligned}
\end{equation}
for all $x_1,x_2 \in \mathbb{Z}$ and
 \begin{equation}\label{441pm522}
 \mathbf{U}(x,x;t)= \mathbf{B}\,\mathbf{U}(x-1,x;t)
\end{equation}
for all $x \in  \mathbb{Z}$. Then, for $x_1=x$ and $x_2 = x+1$, (\ref{748pm51643}) equals
\begin{equation}\label{702pm72}
\begin{aligned}
\frac{d}{dt}\mathbf{U}(x,x+1;t) =~& p\mathbf{U}(x-1,x+1;t) +   p\mathbf{B}\,\mathbf{U}(x-1,x;t)  \\
 &~\,+ q\mathbf{U}(x,x+2;t)  + q\mathbf{B}\,\mathbf{U}(x,x+1;t)- 2\mathbf{U}(x,x+1;t)
\end{aligned}
\end{equation}
which is in form of  (\ref{752pm752}). Thus, if $\mathbf{U}(x_1,x_2;t)$  satisfies  (\ref{748pm51643}) and  the so called, the \textit{boundary condition} (\ref{441pm522}), then it obviously satisfies (\ref{748pm516}) when $x_1< x_2-1$, and satisfies (\ref{752pm752}) when $x_1=x$ and $x_2 = x+1$.

By  the ansatz of the separation of variables $U_{(Y,\nu)}(X,\pi;t)= U_{(Y,\nu)}(X,\pi)T(t)$ and letting $\mathbf{U}(X) = \mathbf{U}(x_1,x_2)$ be the matrix such that $\mathbf{U}(X;t) =\mathbf{U}(X) T(t)$ and then substituting it into (\ref{748pm51643}), we obtain $T(t) = e^{\varepsilon t}$  for some constant  $\varepsilon$ (a constant with respect to $t,x_1,x_2$). The equations for spatial variables of (\ref{748pm51643}) and (\ref{441pm522}) are
\begin{equation}\label{644pm621}
\begin{aligned}
\varepsilon\mathbf{U}(x_1,x_2) = ~& p\mathbf{U}(x_1-1,x_2) + p\mathbf{U}(x_1,x_2-1)  \\
&~~~+ q\mathbf{U}(x_1+1,x_2) + q\mathbf{U}(x_1,x_2+1) - 2\mathbf{U}(x_1,x_2)
\end{aligned}
\end{equation}
 and
 \begin{equation}\label{701pm621}
 \mathbf{U}(x,x)= \mathbf{B}\,\mathbf{U}(x-1,x).
\end{equation}

The Bethe ansatz solution of the matrix equation (\ref{644pm621}) is given by
\begin{equation}\label{442pm522}
\mathbf{U}(x_1,x_2)= \mathbf{A}_{12}\xi_1^{x_1}\xi_2^{x_2} +\mathbf{A}_{21}\xi_2^{x_1}\xi_1^{x_2}
\end{equation}
for any  nonzero complex numbers $\xi_1$ and $\xi_2$. Here, $\mathbf{A}_{12}$ and $\mathbf{A}_{21}$ are $4\times 4$ matrices  similar in structure to $\mathbf{U}(x_1,x_2)$, meaning that the positions of the zero elements in  $\mathbf{A}_{12}$, $\mathbf{A}_{21}$ and $\mathbf{U}(x_1,x_2)$ correspond. The nonzero elements of $\mathbf{A}_{12}$ and $\mathbf{A}_{21}$ are constants (with respect to $t, x_1,x_2$). By substituting (\ref{442pm522}) in (\ref{644pm621}), we  find
\begin{equation}\label{645pm1027}
\varepsilon =\frac{p}{\xi_1} + \frac{p}{\xi_2} +  q{\xi_1} + q{\xi_2}- 2.
\end{equation}
Additionally, since (\ref{442pm522}) must satisfy (\ref{701pm621}), substituting (\ref{442pm522}) into (\ref{701pm621}) yields
\begin{equation}\label{104am1019}
\mathbf{A}_{21} =  -\big(\mathbf{I} - \mathbf{B}\tfrac{1}{\xi_2}\big)^{-1}\big(\mathbf{I} - \mathbf{B}\tfrac{1}{\xi_1}\big)\mathbf{A}_{12}.
\end{equation}
It is important to note that  $\mathbf{I} - \mathbf{B}\tfrac{1}{\xi}$ is  block-diagonal, allowing us to  compute
\begin{equation}\label{352am728}
 -\big(\mathbf{I} - \mathbf{B}\tfrac{1}{\xi_2}\big)^{-1}\big(\mathbf{I} - \mathbf{B}\tfrac{1}{\xi_1}\big)
\end{equation}
\textit{block-wise}.  Specifically, for the $2\times 2$ blocks in the center of (\ref{352am728}),
we have
\begin{equation}\label{142am1019}
-\Bigg[\mathbf{I}_2 - \kbordermatrix{
& 12 & 21\\
12& 0 &0 \\
21& \tfrac{1}{\xi_{2}}& \tfrac{1}{\xi_{2}}
}\Bigg]^{-1}\Bigg[\mathbf{I}_2 - \kbordermatrix{
& 12 & 21\\
12& 0 &0 \\
21& \tfrac{1}{\xi_{1}}& \tfrac{1}{\xi_{1}}
}\Bigg] = \kbordermatrix{
& 12 & 21\\
12& -1 &0 \\
21& T_{21}&  S_{21}
}
\end{equation}
where
\begin{equation}\label{529pm625}
S_{21} =  -\frac{1-1/\xi_{1}}{1-1/\xi_{2}}~~\textrm{and}~~T_{21} = \frac{1/\xi_{1}-1/\xi_{2}}{1-1/\xi_{2}}.
\end{equation}
Hence, we obtain
\begin{equation}\label{230am1019}
\begin{aligned}
\mathbf{A}_{21} =
\kbordermatrix{
&11 &12 &21 &22 \\
                                  11    &  S_{21} & 0 & 0 & 0 \\[6pt]
                                  12    &  0 & -1 & 0 &  0\\[6pt]
                                  21    &  0 &T_{21} & S_{21} & 0 \\[6pt]
                                  22    &  0 & 0 & 0 &  S_{21} \\[6pt]
                                    } \mathbf{A}_{12}
\end{aligned}
\end{equation}

as a condition for (\ref{701pm621}) to be satisfied. Thus, $\mathbf{U}(x_1,x_2)e^{\varepsilon t}$ with  (\ref{442pm522}) and  (\ref{230am1019}) satisfies  (\ref{748pm516}) when $x_1<x_2-1$, and  (\ref{702pm72}) when $x_1=x$ and $x_2 = x+1$. (For now, $\mathbf{U}(x_1,x_2)e^{\varepsilon t}$ is not the transition probability yet because we have not shown that it satisfies the initial condition.)
\begin{remark}
A key idea of defining our two-sided model  as in Definition \ref{836pm72}  comes from the observation that the boundary condition for the one-sided multi-species model with long-range jumps (to the right direction) and the boundary condition for the multi-species TASEP (to the left direction) are the same, as  given by (\ref{701pm621}).
\end{remark}
\section{$N=3$}\label{924pm72}
For $N=3$, $\mathbf{P}_Y(X;t)= \mathbf{P}_Y(x_1,x_2,x_3;t)$ is a $27 \times 27$ matrix whose rows and columns are labelled by $111,112,\dots, 333$.  For $X=(x_1,x_2,x_3)$ with $x_1< x_2-1<x_3-2$, we immediately obtain the master equation in matrix form:
\begin{equation}\label{529pm522}
\begin{aligned}
\frac{d}{dt}\mathbf{P}_{Y}(x_1,x_2,x_3;t)=~& p\mathbf{P}_{Y}(x_1-1,x_2,x_3;t) + p\mathbf{P}_{Y}(x_1,x_2-1,x_3;t) + p\mathbf{P}(x_1,x_2,x_3-1;t)\\
+~& q\mathbf{P}_{Y}(x_1+1,x_2,x_3;t) + q\mathbf{P}_{Y}(x_1,x_2+1,x_3;t) + q\mathbf{P}(x_1,x_2,x_3+1;t)\\[3pt]
&\hspace{2cm} - 3\mathbf{P}_{Y}(x_1,x_2,x_3;t).
\end{aligned}
\end{equation}
\subsection{Master equation for $X=(x_1,x_2,x_3)$ with $x_1=x$, $x_2=x+1$ and $x_2<x_3-1$}\label{334am726}
In this case, we could find  the master equation in matrix form by finding all master equations for the elements of  $\mathbf{P}_Y(X;t)$, but we will use another approach. Considering the positions of particles and the definition of the model, we see that the master equation for $(x,x+1, x_3)$  should be in the form of
\begin{equation}\label{346pm526}
\begin{aligned}
\frac{d}{dt}\mathbf{P}_Y(x,x+1,x_3;t) = ~&  p\mathbf{P}_Y(x-1,x+1,x_3;t) + \mathbf{C}\, \mathbf{P}_Y(x-1,x,x_3;t)\\
+~& p\mathbf{P}_Y(x,x+1,x_3-1;t)-3p\mathbf{P}_Y(x,x+1,x_3;t)\\[3pt]
+~&\mathbf{C}_1\mathbf{P}_Y(x,x+1,x_3;t) + q\mathbf{P}_Y(x,x+2,x_3;t)\\[3pt]
+~& q\mathbf{P}_Y(x,x+1,x_3+1;t)-\mathbf{C}_2\mathbf{P}_Y(x,x+1,x_3;t)
\end{aligned}
\end{equation}
for some $27 \times 27$ matrices $\mathbf{C}, \mathbf{C}_1, \mathbf{C}_2$.
\subsubsection{Matrix $\mathbf{C}$}\label{C1}
The  $(ijk,lmn)$-element of $\mathbf{C}$ is the instantaneous rate of transition from state $(x-1,x,x_3,lmn)$ to state $(x,x+1,x_3,ijk)$. In fact, if the transition
\begin{equation}\label{623pm545}
(x-1,x,x_3,lmn)~~\longrightarrow~~(x,x+1,x_3,ijk)
\end{equation}
is permissible, then it occurs with rate $p$ by the definition of the model, and if (\ref{623pm545}) is not permissible, its rate is obviously 0. The permissibility of (\ref{623pm545}) depends on $ijk$ and $lmn$. Since $x_3$ is kept the same during the transition, it means that $k=n$. If $k\neq n$,  (\ref{623pm545}) is not possible so its rate is 0. In the case that $k=n$, we may say that the permissibility of  (\ref{623pm545}) depends on the permissibility of the transition from $(x-1,x,lm)$ to $(x,x+1,ij)$, which is actually  encoded in  (\ref{Bmatrix}). Now, let us view $\mathbf{C}$ as a $9 \times 9$ partitioned matrix whose rows and columns are labelled by $11\cdot, 12\cdot, 13\cdot, 21\cdot, 22\cdot, 23\cdot, 31\cdot, 32\cdot, 33\cdot$, and the $(ij\cdot, lm\cdot)$-element is a $3\times 3$ matrix whose rows are labelled by $ij1,ij2,ij3$ and columns are labelled by $lm1,lm2,lm3$. According to the permissibility of the transition from $(x-1,x,lm)$ to $(x,x+1,ij)$  available in (\ref{Bmatrix}) and the requirement $k=n$ for  (\ref{623pm545}), we obtain
\begin{equation*}
\mathbf{C}= \kbordermatrix{
       & 11\cdot & 12\cdot & 13\cdot & 21\cdot & 22\cdot & 23\cdot & 31\cdot & 32\cdot & 33\cdot\\
     11\cdot & p\mathbf{I} & \mathbf{0} & \mathbf{0} & \mathbf{0} & \mathbf{0} & \mathbf{0} & \mathbf{0} & \mathbf{0}& \mathbf{0}\\
     12\cdot  & \mathbf{0} & \mathbf{0} &  \mathbf{0} &  \mathbf{0} &  \mathbf{0} &  \mathbf{0} & \mathbf{0}  &  \mathbf{0}& \mathbf{0}\\
     13\cdot  & \mathbf{0} &  \mathbf{0} &  \mathbf{0} &  \mathbf{0} &  \mathbf{0} &  \mathbf{0} &  \mathbf{0} & \mathbf{0} & \mathbf{0}\\
      21\cdot & \mathbf{0} &  p\mathbf{I} &  \mathbf{0} &  p\mathbf{I} &  \mathbf{0} &  \mathbf{0} &  \mathbf{0} &  \mathbf{0}& \mathbf{0}\\
     22\cdot  & \mathbf{0} &  \mathbf{0} &  \mathbf{0} &  \mathbf{0} &  p\mathbf{I} &  \mathbf{0} &  \mathbf{0} &  \mathbf{0}& \mathbf{0}\\
     23\cdot & \mathbf{0} &  \mathbf{0} &  \mathbf{0} &  \mathbf{0} &  \mathbf{0} &  \mathbf{0} &  \mathbf{0} &  \mathbf{0}& \mathbf{0}\\
     31\cdot  & \mathbf{0} &  \mathbf{0} &  p\mathbf{I} &  \mathbf{0} &  \mathbf{0} &  \mathbf{0} &  p\mathbf{I} & \mathbf{0} & \mathbf{0}\\
    32\cdot   & \mathbf{0} &  \mathbf{0} &  \mathbf{0} &  \mathbf{0} &  \mathbf{0} &  p\mathbf{I} &  \mathbf{0} &  p\mathbf{I}& \mathbf{0}\\
     33\cdot & \mathbf{0} &  \mathbf{0} &  \mathbf{0} &  \mathbf{0} &  \mathbf{0} &  \mathbf{0} &  \mathbf{0} &  \mathbf{0}& p\mathbf{I}\\
 }.
\end{equation*}
If we let
\begin{equation}\label{646pm623}
\mathbf{B}= \kbordermatrix{
       & 11 & 12 & 13 & 21 & 22& 23 & 31 & 32 & 33\\
     11 & 1 &0 & 0 & 0 & 0 & 0 & 0 & 0& 0\\
     12  & 0 & 0 &  0 &  0 &  0 &  0 &0  &  0& 0\\
     13  & 0 & 0 &  0 & 0&  0 &  0 & 0 & 0 & 0\\
      21&0 &  1 &  0 &  1 & 0 &  0 &  0 & 0& 0\\
     22  & 0 &  0 &  0 &  0 &  1 &  0 &  0 & 0& 0\\
     23 & 0 &  0 &  0 &  0& 0 &  0 &  0 &  0& 0\\
     31 & 0 & 0 & 1&  0 &  0 & 0 & 1 & 0 & 0\\
    32   & 0 &  0 & 0 &  0 & 0 &  1 &  0 & 1&0\\
     33 & 0 &  0 & 0 & 0 &  0 &  0 & 0 &  0& 1\\
 },
\end{equation}
then we have a succinct expression:
\begin{equation}\label{1059am73}
\mathbf{C}= p\mathbf{B} \otimes \mathbf{I}.
\end{equation}
\begin{remark}
We have already used the notation $\mathbf{B}$ in (\ref{Bmatrix}) for $N=2$, but here for $N=3$ we use it again because later the notation $\mathbf{B}$ will be defined for general $N$ consistently with (\ref{Bmatrix}) and  (\ref{646pm623}).
\end{remark}
\subsubsection{Matrices $\mathbf{C}_1$ and $\mathbf{C}_2$}\label{941pm725}
The  $(ijk,lmn)$-element of $\mathbf{C}_1$ is the instantaneous rate of transition
\begin{equation}\label{623pm5458}
(x,x+1,x_3,lmn)~~\longrightarrow~~(x,x+1,x_3,ijk)
\end{equation}
which is due to the interchange of positions of particles at $x$ and $x+1$. Hence, the instantaneous rate of (\ref{623pm5458}) is $q$ if the transition is permissible, and otherwise, the rate is 0. Since $x+1< x_3-1$ in (\ref{623pm5458}), it means that $k=n$, and we see that the permissibility of (\ref{623pm5458}) depends on the permissibility of
\begin{equation*}
(x,x+1,lm)~~\longrightarrow~~(x,x+1,ij)
\end{equation*}
which is actually encoded in  (\ref{Bmatrix2}). Using a similar argument as in Section \ref{C1}, if we let
\begin{equation}\label{646pm6236}
\mathbf{B}_1= \kbordermatrix{
       & 11 & 12 & 13 & 21 & 22& 23 & 31 & 32 & 33\\
     11 & 0 &0 & 0 & 0 & 0 & 0 & 0 & 0& 0\\
     12  & 0 & 0 &  0 &  0 &  0 &  0 &0  &  0& 0\\
     13  & 0 & 0 &  0 & 0&  0 &  0 & 0 & 0 & 0\\
      21&0 &  1 &  0 &  0 & 0 &  0 &  0 & 0& 0\\
     22  & 0 &  0 &  0 &  0 &  0 &  0 &  0 & 0& 0\\
     23 & 0 &  0 &  0 &  0& 0 &  0 &  0 &  0& 0\\
     31 & 0 & 0 & 1&  0 &  0 & 0 & 0 & 0 & 0\\
    32   & 0 &  0 & 0 &  0 & 0 &  1 &  0 & 0&0\\
     33 & 0 &  0 & 0 & 0 &  0 &  0 & 0 &  0& 0\\
 },
\end{equation}
 we have
\begin{equation*}
\mathbf{C}_1 = q\mathbf{B}_1 \otimes \mathbf{I}.
\end{equation*}
\begin{remark}
We have already used the notation $\mathbf{B}_1$ in (\ref{Bmatrix2}) for $N=2$, but here for $N=3$ we use it again because later the notation $\mathbf{B}_1$ will be defined for general $N$ consistently with (\ref{Bmatrix2}) and  (\ref{646pm6236}).
\end{remark}
The matrix $\mathbf{C}_2$ is a diagonal matrix such that the $(ijk,ijk)$-element is the instantaneous rate of jumping out of the state $(x,x+1,x_3,ijk)$ due to particles' jumps to the left. The rate of the jump of the particle $k$ at $x_3$ to the left direction is $q$ and the rate of the jump of $i$ at $x$ and $j$ at $x+1$ is encoded in (\ref{Bmatrix3}). Hence, if we let
\begin{equation}\label{646pm62361}
\mathbf{B}_2= \kbordermatrix{
       & 11 & 12 & 13 & 21 & 22& 23 & 31 & 32 & 33\\
     11 & 1 &0 & 0 & 0 & 0 & 0 & 0 & 0& 0\\
     12  & 0 & 2 &  0 &  0 &  0 &  0 &0  &  0& 0\\
     13  & 0 & 0 &  2 & 0&  0 &  0 & 0 & 0 & 0\\
      21&0 &  0 &  0 &  1 & 0 &  0 &  0 & 0& 0\\
     22  & 0 &  0 &  0 &  0 &  1 &  0 &  0 & 0& 0\\
     23 & 0 &  0 &  0 &  0& 0 &  2 &  0 &  0& 0\\
     31 & 0 & 0 & 0&  0 &  0 & 0 & 1 & 0 & 0\\
    32   & 0 &  0 & 0 &  0 & 0 &  0 &  0 & 1&0\\
     33 & 0 &  0 & 0 & 0 &  0 &  0 & 0 &  0& 1\\
 },
\end{equation}
then we obtain
\begin{equation*}
\mathbf{C}_2 = q\mathbf{B}_2 \otimes \mathbf{I} + q\mathbf{I}_{27}.
\end{equation*}
\begin{remark}
We have already used the notation $\mathbf{B}_2$ in (\ref{Bmatrix3}) for $N=2$, but here for $N=3$ we use it again because later the notation $\mathbf{B}_2$ will be defined for general $N$ consistently with (\ref{Bmatrix3}) and  (\ref{646pm62361}).
\end{remark}
Moreover, observing that $\mathbf{B}_2 + \mathbf{B} - \mathbf{B}_1 = 2\mathbf{I}_9$,  we see that (\ref{346pm526}) is equivalent to
\begin{equation}\label{346pm52623}
\begin{aligned}
&\frac{d}{dt}\mathbf{P}_Y(x,x+1,x_3;t)  =  p\mathbf{P}_Y(x-1,x+1,x_3;t) + p(\mathbf{B} \otimes \mathbf{I})\, \mathbf{P}_Y(x-1,x,x_3;t) \\
&\hspace{1cm} + p\mathbf{P}_Y(x,x+1,x_3-1;t) +q(\mathbf{B} \otimes \mathbf{I})\mathbf{P}_Y(x,x+1,x_3;t) + q\mathbf{P}_Y(x,x+2,x_3;t)\\[3pt]
&\hspace{1cm} +q\mathbf{P}_Y(x,x+1,x_3+1;t)-3\mathbf{P}_Y(x,x+1,x_3;t).
\end{aligned}
\end{equation}
\subsection{Master equation for $X=(x_1,x_2,x_3)$ with $x_1+1<x_2$, $x_2=x$ and $x_3=x+1$}\label{453am726}
Considering the positions of particles and the definition of the model, we see that the master equation  should be in the form:
\begin{equation}\label{911pm1016}
\begin{aligned}
\frac{d}{dt}\mathbf{P}_Y(x_1,x,x+1;t) = ~&  p\mathbf{P}_Y(x_1-1,x,x+1;t) + p\mathbf{P}_Y(x_1,x-1,x+1;t)\\
+~& \mathbf{E}\, \mathbf{P}_Y(x_1,x-1,x;t)-3p\mathbf{P}_Y(x_1,x,x+1;t)\\[3pt]
+~&  q\mathbf{P}_Y(x_1+1,x,x+1;t)+ \mathbf{E}_1\mathbf{P}_Y(x_1,x,x+1;t)\\[1pt]
+~& q\mathbf{P}_Y(x_1,x,x+2;t)-\mathbf{E}_2\mathbf{P}_Y(x_1,x,x+1;t)
\end{aligned}
\end{equation}
for some $27 \times 27$ matrices $\mathbf{E}, \mathbf{E}_1, \mathbf{E}_2$.

First, we find the matrix $\mathbf{E}$. The $(ijk,lmn)$-element of $\mathbf{E}$ is the instantaneous rate of  the transition
\begin{equation}\label{623pm5456}
(x_1,x-1,x,lmn)~~\longrightarrow~~(x_1,x,x+1,ijk).
\end{equation}
If the transition (\ref{623pm5456}) is permissible, then it occurs with rate $p$ by the definition of the model. If (\ref{623pm5456}) is not permissible, its rate is obviously $0$.     Since $x_1$ is kept the same during the transition (\ref{623pm5456}), it means that $i=l$. If  $i\neq l$, then  (\ref{623pm5456}) is not possible, so its rate is  0. In the case that $i=l$,  we may say that the permissibility of   (\ref{623pm5456}) depends on the permissibility of the transition from $(x-1,x,mn)$ to $(x,x+1,jk)$, which is actually encoded in the matrix $\mathbf{B}$ in (\ref{646pm623}).

Now, let us view $\mathbf{E}$ as a $3\times 3$ partitioned matrix whose rows and columns are labelled by $1**, 2**, 3\**$. The $(i**, l**)$-element is a $9 \times 9$ matrix whose rows and columns are labelled by $i11,i12,\dots, i33$ and $l11,l12,\dots, l33$, respectively. Then, using $\mathbf{B}$ in (\ref{646pm623}) and the requirement $i =l$ for  the transition (\ref{623pm5456}), we obtain
\begin{equation}\label{1100am73}
\mathbf{E}=\kbordermatrix{
       & 1** & 2** & 3** \\
    1** & p\mathbf{B}  & \mathbf{0} & \mathbf{0}\\
    2**  & \mathbf{0} & p\mathbf{B}  & \mathbf{0}  \\
    3**  & \mathbf{0} & \mathbf{0} &  p\mathbf{B} \\
 } = \mathbf{I}\otimes p\mathbf{B}.
 \end{equation}
Using similar arguments as in Section \ref{941pm725}, we obtain:
 \begin{equation*}
 \mathbf{E}_1 =  \mathbf{I} \otimes q\mathbf{B}_1
 \end{equation*}
 with  $\mathbf{B}_1$ in (\ref{646pm6236}), and
 \begin{equation*}
\mathbf{E}_2 = \mathbf{I} \otimes q\mathbf{B}_2  + q\mathbf{I}_{27}
\end{equation*}
with $\mathbf{B}_2$ in (\ref{646pm62361}). Moreover, observing that $\mathbf{B}_2 + \mathbf{B} - \mathbf{B}_1 = 2\mathbf{I}_9$, we see that (\ref{911pm1016}) is equivalent to
\begin{equation}\label{346pm526232}
\begin{aligned}
&\frac{d}{dt}\mathbf{P}_Y(x_1,x,x+1;t) =   p\mathbf{P}_Y(x_1-1,x,x+1;t) + p\mathbf{P}_Y(x_1,x-1,x+1;t)\\
&\hspace{1cm}+ p(\mathbf{I} \otimes \mathbf{B}) \, \mathbf{P}_Y(x_1,x-1,x;t)+  q\mathbf{P}_Y(x_1+1,x,x+1;t)\\[3pt]
&\hspace{1cm} +q(\mathbf{I} \otimes \mathbf{B})\mathbf{P}_Y(x_1,x,x+1;t)+ q\mathbf{P}_Y(x_1,x,x+2;t)-3\mathbf{P}_Y(x_1,x,x+1;t).
\end{aligned}
\end{equation}
\subsection{Master equation for $X=(x,x+1,x+2)$}\label{802pm728}
Considering all possible states that can transition to states with $(x,x+1,x+2)$, we see that  the master equation should be in the form:
\begin{equation}\label{350pm526}
\begin{aligned}
\frac{d}{dt}\mathbf{P}_Y(x,x+1,x+2;t) =~& p\mathbf{P}_Y(x-1,x+1,x+2;t) + \mathbf{A} \mathbf{P}_Y(x-1,x,x+2;t) \\
+~&\mathbf{A}_1\,\mathbf{P}_Y(x-1,x,x+1;t) -3p\,\mathbf{P}_Y(x,x+1,x+2;t) \\
+~& q\mathbf{P}_Y(x,x+1,x+3;t) + \mathbf{A}_2 \mathbf{P}_Y(x,x+1,x+2;t)\\
-~& \mathbf{A}_3 \mathbf{P}_Y(x,x+1,x+2;t)
\end{aligned}
\end{equation}
for some  $27 \times 27$ matrices  $\mathbf{A},\mathbf{A}_1, \mathbf{A}_2, \mathbf{A}_3$.

The $(ijk, lmn)$-element of $\mathbf{A}$ is the instantaneous rate of transition from  $(x-1,x,x+2,lmn)$ to  $(x,x+1,x+2,ijk)$.  Using the same argument as in Section \ref{941pm725}, we find that $\mathbf{A} = \mathbf{C} =  p\mathbf{B} \otimes \mathbf{I}$.

The $(ijk, lmn)$-element of $\mathbf{A}_1$ is the instantaneous rate of transition
\begin{equation}\label{627pm624}
(x-1,x,x+1,lmn) ~~ \longrightarrow~~(x,x+1,x+2,ijk)
\end{equation}
We could find all matrix elements of  $\mathbf{A}_1$ by checking the permissibility of (\ref{627pm624}), but we will take a different approach by considering an alternative perspective on the definition of our model. Specifically, we will interpret the rule of the long-range jump as consecutive two-body interactions.

\begin{figure}[h]
  \centering
  \begin{tikzpicture}
    % circles with numbers inside
    \node[circle, draw, minimum size=0.35cm, inner sep = 0] (left1) at (0,0) {$1$};
    \node[circle, draw, minimum size=0.35cm, inner sep = 0] (left2) at (1,0) {$2$};

     % curved arrow
    \draw[->, thick, bend left=-45] ([yshift=0.2cm]left2.north) to ([yshift=0.2cm]left1.north);

    % Arrow pointing to the right
    \draw[->] (2.3,0) -- (2.7,0);
    \draw[->] (6.3,0) -- (6.7,0);

    % circles with numbers inside
    \node[circle, draw, minimum size=0.35cm, inner sep = 0] at (3.82,0) {$2$};
    \node[circle, draw, minimum size=0.35cm, inner sep = 0] at (4.18,0) {$1$};
    \node[circle, draw, dashed, thick, minimum size=0.35cm, inner sep = 0] at (5,0) {};

    % circles with numbers inside
    \node[circle, draw, minimum size=0.35cm, inner sep = 0] at (8,0) {$2$};
    \node[circle, draw, minimum size=0.35cm, inner sep = 0] at (9,0) {$1$};

  \end{tikzpicture}

  \caption{Transition from $(x,x+1,12)$ to $(x,x+1,21)$}
  \label{135pm726}
\end{figure}
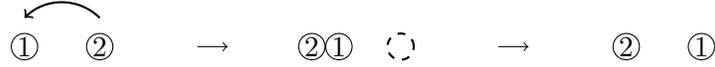

\begin{definition}\label{636pm624} (equivalent to Definition \ref{836pm72})   A particle  $l$ at position $x$ waits an exponential time with rate 1 and then jumps to $x+1$ with probability $p$ or to $x-1$ with probability $q=1-p$. When the particle jumps to $x \pm 1$, if $x \pm 1$ is already occupied by a particle $l'$, then $x \pm 1$ instantaneously accommodates both $l$ and $l'$  such that $\max(l,l')$  is positioned to the left of $\min(l,l')$ at $x \pm 1$, and $\min(l,l')$ immediately jumps to $(x \pm 1) +1$.
\end{definition}
Figure \ref{130pm726} and \ref{135pm726} illustrate some transitions as described in Definition \ref{636pm624}.
\begin{figure}[h]
  \centering
  \begin{tikzpicture}
    % circles with numbers inside
    \node[circle, draw, minimum size=0.35cm, inner sep = 0] (left1) at (0,0) {$1$};
    \node[circle, draw, minimum size=0.35cm, inner sep = 0] (left2) at (1,0) {$2$};
    \node[circle, draw, minimum size=0.35cm, inner sep = 0] (left3) at (2,0) {$3$};

     % curved arrow
    \draw[->, thick, bend left=45] ([yshift=0.2cm]left1.north) to ([yshift=0.2cm]left2.north);

    % Arrow pointing to the right
    \draw[->] (2.8,0) -- (3.2,0);
    \draw[->] (6.8,0) -- (7.2,0);
    \draw[->] (10.8,0) -- (11.2,0);

    % circles with numbers inside
    \node[circle, draw, dashed, thick, minimum size=0.35cm, inner sep = 0] at (4,0) {};
    \node[circle, draw, minimum size=0.35cm, inner sep = 0] at (4.82,0) {$2$};
    \node[circle, draw, minimum size=0.35cm, inner sep = 0] at (5.18,0) {$1$};
    \node[circle, draw, minimum size=0.35cm, inner sep = 0] at (6,0) {$3$};

    % circles with numbers inside
    \node[circle, draw, dashed, thick, minimum size=0.35cm, inner sep = 0] at (8,0) {};
    \node[circle, draw, minimum size=0.35cm, inner sep = 0] at (9,0) {$2$};
    \node[circle, draw, minimum size=0.35cm, inner sep = 0] at (9.82,0) {$3$};
    \node[circle, draw, minimum size=0.35cm, inner sep = 0] at (10.18,0) {$1$};

    % circles with numbers inside
    \node[circle, draw, dashed, thick, minimum size=0.35cm, inner sep = 0] at (12,0) {};
    \node[circle, draw, minimum size=0.35cm, inner sep = 0] at (13,0) {$2$};
    \node[circle, draw, minimum size=0.35cm, inner sep = 0] at (14,0) {$3$};
    \node[circle, draw, minimum size=0.35cm, inner sep = 0] at (15,0) {$1$};

  \end{tikzpicture}

  \caption{Transition from $(x,x+1,x+2,123)$ to $(x+1,x+2,x+3,231)$}
  \label{130pm726}
\end{figure}
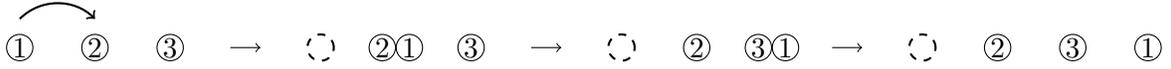

According to Definition \ref{636pm624}, the transition (\ref{627pm624}) can be interpreted as
\begin{equation}\label{628pm624}
\begin{aligned}
&(x-1,x,x+1,lmn) ~ \xrightarrow{~~a~~}~(x,x,x+1,l'm'n')\\[5pt]
&\hspace{1.5cm}~ \xrightarrow{~~b~~}~(x,x+1,x+1,l''m''n'')~ \xrightarrow{~~c~~}~(x,x+1,x+2,ijk)
\end{aligned}
\end{equation}
where $\xrightarrow{~~b~~}$ and $\xrightarrow{~~c~~}$ occur instantaneously. We see that  for (\ref{628pm624}) to be possible it must be that $n = n', l' = l''$ and $l''m''n''=ijk$. Using the arguments in Section \ref{334am726} and \ref{453am726}, the permissibility of $\xrightarrow{~~a~~}$  is represented by $\mathbf{B} \otimes \mathbf{I}$ and  the permissibility of $\xrightarrow{~~b~~}$  is represented by $\mathbf{I}\otimes \mathbf{B}$. It is obvious that $\xrightarrow{~~c~~}$ is represented by the identity matrix.  Hence, we conclude that L
\begin{equation}\label{628pm73}
\mathbf{A}_1 =  p(\mathbf{I}\otimes \mathbf{B})(\mathbf{B} \otimes \mathbf{I}).
\end{equation}
The $(ijk,lmn)$-element of $\mathbf{A}_2$ is the instantaneous rate of transition from  $(x,x+1,x+2,lmn)$ to  $(x,x+1,x+2,ijk)$. This transition is possible when either $m>l$ or $n>m$ according to the rule of the mTASEP to the left direction. Hence, using the arguments in Section \ref{941pm725}, we have
\begin{equation*}
\mathbf{A}_2 = q(\mathbf{B}_1 \otimes \mathbf{I}) +  q(\mathbf{I}\otimes \mathbf{B}_1).
\end{equation*}
$\mathbf{A}_3$ is a diagonal matrix such that the $(ijk,ijk)$-element represents the rate of jumping out of the state $(x,x+1,x+2,ijk)$ due to particles' jumps to the left direction. Recall that $q\mathbf{B}_2 \otimes \mathbf{I}$ describes the rate of jumping out of the state $(x,x+1,x_3,ijk)$ with $x+1 < x_3-1$ due to the first particle's jump  to the left  and the second particle's jump  to the left (which is possible only when $j>i$) and $\mathbf{I} \otimes q\mathbf{B}_2$  describes the rate of jumping out of the state $(x_1,x,x+1,ijk)$ with $x_1<x-1$ due to the second particle's jump to the left  and third particle's jump to the left (which is possible only when $k>j$). Hence, it must be that
\begin{equation*}
\mathbf{A}_3 = q\mathbf{B}_2 \otimes \mathbf{I}  + \mathbf{I} \otimes q\mathbf{B}_2 - q\mathbf{I}_{27}.
\end{equation*}
Again, using that $\mathbf{B}_2 + \mathbf{B} - \mathbf{B}_1 = 2\mathbf{I}_{9}$, we see that (\ref{350pm526}) is equivalent to
\begin{equation}\label{350pm52632}
\begin{aligned}
&\frac{d}{dt}\mathbf{P}_Y(x,x+1,x+2;t) = p\mathbf{P}_Y(x-1,x+1,x+2;t) + p(\mathbf{B} \otimes \mathbf{I}) \mathbf{P}_Y(x-1,x,x+2;t) \\
&\hspace{1cm}+p(\mathbf{I}\otimes \mathbf{B})(\mathbf{B} \otimes \mathbf{I})\,\mathbf{P}_Y(x-1,x,x+1;t) +q\mathbf{P}_Y(x,x+1,x+3;t) \\[4pt]
&\hspace{1cm}+  q(\mathbf{B} \otimes \mathbf{I} +  \mathbf{I}\otimes \mathbf{B}) \mathbf{P}_Y(x,x+1,x+2;t)- 3 \mathbf{P}_Y(x,x+1,x+2;t).
\end{aligned}
\end{equation}
\subsection{Two-body reducibility}
We obtained four different forms  of the master equation  (\ref{529pm522}), (\ref{346pm52623}), (\ref{346pm526232}), and (\ref{350pm52632}) for $\mathbf{P}_Y(X;t)$, depending on $X=(x_1,x_2,x_3)$. Now, we will show that (\ref{346pm52623}), (\ref{346pm526232}), and (\ref{350pm52632}) can be derived from the equation in the form of  (\ref{529pm522}) and the \textit{boundary conditions} that describe two-particle interactions.

We naturally extend the notations $\mathbf{U}(X;t)$ and  $U_{(Y,\nu)}(X,\pi;t)$ for $N=2$ in Section \ref{646pm624}  to $N=3$. Suppose that $\mathbf{U}(X;t) = \mathbf{U}(x_1,x_2,x_3;t)$ satisfies
 \begin{equation}\label{209pm10171}
\begin{aligned}
\frac{d}{dt}\mathbf{U}(x_1,x_2,x_3;t)=~& p\mathbf{U}(x_1-1,x_2,x_3;t) + p\mathbf{U}(x_1,x_2-1,x_3;t) + p\mathbf{U}(x_1,x_2,x_3-1;t)\\
+~&q\mathbf{U}(x_1+1,x_2,x_3;t) + q\mathbf{U}(x_1,x_2+1,x_3;t) + q\mathbf{U}(x_1,x_2,x_3+1;t)\\[3pt]
&\hspace{2cm} - 3\mathbf{U}(x_1,x_2,x_3;t)
\end{aligned}
\end{equation}
for all $x_1,x_2,x_3 \in \mathbb{Z}$. If $\mathbf{U}(x_1,x_2,x_3;t)$ also satisfies
\begin{equation}\label{1250pm625}
(\mathbf{B}\otimes  \mathbf{I}) \mathbf{U}(x-1,x,x_3;t) =  \mathbf{U}(x,x,x_3;t)~~~\textrm{for all $x,x_3 \in \mathbb{Z}$,}
\end{equation}
then for $x_1=x, x_2 = x+1, x_2<x_3-1$, (\ref{209pm10171}) becomes
  \begin{equation}\label{209pm101711}
\begin{aligned}
&\frac{d}{dt}\mathbf{U}(x,x+1,x_3;t)  =  p\mathbf{U}(x-1,x+1,x_3;t) + p(\mathbf{B} \otimes \mathbf{I})\, \mathbf{U}(x-1,x,x_3;t) \\
&\hspace{1cm} + p\mathbf{U}(x,x+1,x_3-1;t) +q(\mathbf{B} \otimes \mathbf{I})\mathbf{U}(x,x+1,x_3;t) + q\mathbf{U}(x,x+2,x_3;t)\\[3pt]
&\hspace{1cm} +q\mathbf{U}(x,x+1,x_3+1;t)-3\mathbf{U}(x,x+1,x_3;t)
\end{aligned}
\end{equation}
which is in the form of (\ref{346pm52623}). Similarly, if  $\mathbf{U}(x_1,x_2,x_3;t)$ satisfies
\begin{equation}\label{1251pm625}
( \mathbf{I}\otimes \mathbf{B})\mathbf{U}(x_1,x-1,x;t) =\mathbf{U}(x_1,x,x;t)~~~\textrm{for all $x_1,x \in \mathbb{Z},$}
\end{equation}
then (\ref{209pm10171}) becomes
 \begin{equation}\label{209pm1017111}
\begin{aligned}
&\frac{d}{dt}\mathbf{U}(x_1,x,x+1;t) =   p\mathbf{U}(x_1-1,x,x+1;t) + p\mathbf{U}(x_1,x-1,x+1;t)\\
&\hspace{1cm}+ p(\mathbf{I} \otimes \mathbf{B}) \, \mathbf{U}(x_1,x-1,x;t)+  q\mathbf{U}(x_1+1,x,x+1;t)\\[3pt]
&\hspace{1cm} +q(\mathbf{I} \otimes \mathbf{B})\mathbf{U}(x_1,x,x+1;t)+ q\mathbf{U}(x_1,x,x+2;t)-3\mathbf{U}(x_1,x,x+1;t).
\end{aligned}
\end{equation}
which is in the form of (\ref{346pm526232}).

Now, we find a condition for (\ref{209pm10171}) to be in the form of (\ref{350pm52632}). For this purpose, we substitute $x_1=x,x_2 = x+1, x_3 = x+2$ in (\ref{209pm10171}) and compare it with (\ref{350pm52632}). Then, we can see that
\begin{equation}\label{914pm1016}
\begin{aligned}
& p\mathbf{U}(x,x,x+2;t) + p\mathbf{U}(x,x+1,x+1;t) +  q\mathbf{U}(x+1,x+1,x+2;t)  + q\mathbf{U}(x,x+2,x+2;t) \\[5pt]
&\hspace{0.5cm}=\,p(\mathbf{B}\otimes  \mathbf{I}) \mathbf{U}(x-1,x,x+2;t) +p(\mathbf{I}\otimes \mathbf{B})(\mathbf{B} \otimes \mathbf{I})\,\mathbf{U}(x-1,x,x+1;t) \\[5pt]
&\hspace{0.5cm}+\,q(\mathbf{B} \otimes \mathbf{I} +  \mathbf{I}\otimes \mathbf{B})\mathbf{U}(x,x+1,x+2;t)
\end{aligned}
\end{equation}
for all $x \in \mathbb{Z}$ if and only if  (\ref{209pm10171}) is in the form of (\ref{350pm52632}). To satisfy (\ref{914pm1016}), we set:
\begin{equation}\label{549pm726}
\begin{aligned}
& p\mathbf{U}(x,x,x+2;t) + p\mathbf{U}(x,x+1,x+1;t) = \\
&\hspace{2cm}p(\mathbf{B}\otimes  \mathbf{I}) \mathbf{U}(x-1,x,x+2;t) +p(\mathbf{I}\otimes \mathbf{B})(\mathbf{B} \otimes \mathbf{I})\,\mathbf{U}(x-1,x,x+1;t)
\end{aligned}
\end{equation}
and
\begin{equation}\label{550pm726}
\begin{aligned}
 &q\mathbf{U}(x+1,x+1,x+2;t)  + q\mathbf{U}(x,x+2,x+2;t) = \\
 &\hspace{4cm} q(\mathbf{B} \otimes \mathbf{I} +  \mathbf{I}\otimes \mathbf{B})\mathbf{U}(x,x+1,x+2;t).
\end{aligned}
\end{equation}
Now, we will show that  (\ref{549pm726}) and (\ref{550pm726})  are satisfied if the \textit{boundary conditions} (\ref{1250pm625}) and (\ref{1251pm625})  are satisfied. It is easy to see that (\ref{550pm726}) is satisfied by  (\ref{1250pm625}) and (\ref{1251pm625}). Let us show that  (\ref{549pm726}) is satisfied by (\ref{1250pm625}) and (\ref{1251pm625}). The first term on the left side of (\ref{549pm726}) is the same as the first term on the right side of (\ref{549pm726}) by (\ref{1250pm625}). Also, the second term on the left side of (\ref{549pm726}) is the same as the second term on the right side of (\ref{549pm726}) because
\begin{equation*}
( \mathbf{I}\otimes \mathbf{B})(\mathbf{B}\otimes  \mathbf{I})\mathbf{U}(x-1,x,x+1;t) = (\mathbf{I}\otimes \mathbf{B})\mathbf{U}(x,x,x+1;t)
\end{equation*}
by (\ref{1250pm625}), and
\begin{equation*}
(\mathbf{I}\otimes \mathbf{B})\mathbf{U}(x,x,x+1;t) = \mathbf{U}(x,x+1,x+1;t)
\end{equation*}
by (\ref{1251pm625}).
\subsection{Bethe ansatz solution}
In the previous section, we showed that if $\mathbf{U}(X;t)$ satisfies (\ref{209pm10171}), (\ref{1250pm625}) and (\ref{1251pm625}), then it satisfies the master equation  (\ref{529pm522}), (\ref{346pm52623}), (\ref{346pm526232}), or (\ref{350pm52632}) for each $X$ in the corresponding physical region. In this section, we find such a solution $\mathbf{U}(X;t)$ that satisfies (\ref{209pm10171}), (\ref{1250pm625}) and (\ref{1251pm625}). Extending the idea and the notations in Section \ref{646pm624},  we separate variables in the form $\textbf{U}(X;t) = \textbf{U}(X)T(t)$. This leads to the following spatial equations derived from (\ref{209pm10171}), (\ref{1250pm625}) and (\ref{1251pm625}):
\begin{equation}\label{644pm6218}
\begin{aligned}
\varepsilon\mathbf{U}(x_1,x_2,x_3) =~& p\mathbf{U}(x_1-1,x_2,x_3) + p\mathbf{U}(x_1,x_2-1,x_3) + p\mathbf{U}(x_1,x_2,x_3-1) \\
+~& q\mathbf{U}(x_1+1,x_2,x_3) + q\mathbf{U}(x_1,x_2+1,x_3) + q\mathbf{U}(x_1,x_2,x_3+1)\\
-~& 3\mathbf{U}(x_1,x_2,x_3)
\end{aligned}
\end{equation}
and
\begin{eqnarray}
(\mathbf{B}\otimes  \mathbf{I}) \mathbf{U}(x-1,x,x_3)& = & \mathbf{U}(x,x,x_3),\label{408pm625} \\
( \mathbf{I}\otimes \mathbf{B})\mathbf{U}(x_1,x-1,x)& = &\mathbf{U}(x_1,x,x). \label{409pm625}
\end{eqnarray}
The Bethe ansatz solution for  (\ref{644pm6218}) is given by
\begin{equation}\label{1118-225}
\mathbf{U}(x_1,x_2,x_3) = \sum_{\sigma\in S_3}\mathbf{A}_{\sigma}\xi_{\sigma(1)}^{x_1}\xi_{\sigma(2)}^{x_2}\xi_{\sigma(3)}^{x_3}
\end{equation}
where $\xi_1,\xi_2,\xi_3$ are nonzero complex numbers. The sum is taken over all permutations $\sigma=\sigma(1)\sigma(2)\sigma(3)$ in the symmetric group $S_3$, and $\mathbf{A}_{\sigma}$ is a $27 \times 27$ matrix with the same structure as  $\mathbf{U}(x_1,x_2,x_3)$ (that is, zero elements in  $\mathbf{A}_{\sigma}$ and  $\mathbf{U}(x_1,x_2,x_3)$ occupy the same positions), whose nonzero elements are constants with respect to $t,x_1,x_2,x_3$. Substituting (\ref{1118-225}) into (\ref{644pm6218}), we find that
\begin{equation*}
\varepsilon = \frac{p}{\xi_1} + \frac{p}{\xi_2} + \frac{p}{\xi_3} + q\xi_1 +q\xi_2 + q\xi_3 - 3.
\end{equation*}

We now substitute (\ref{1118-225}) into (\ref{408pm625}) and (\ref{409pm625}) to derive conditions on $\mathbf{A}_{\sigma}$ that ensure (\ref{408pm625}) and (\ref{409pm625}) are satisfied. First, substituting (\ref{1118-225}) into (\ref{408pm625}) and simplifying,  we obtain
\begin{equation}\label{823pm61}
\begin{aligned}
&(\mathbf{A}_{123}+\mathbf{A}_{213})+ (\mathbf{A}_{132}+\mathbf{A}_{312})+ (\mathbf{A}_{231}+\mathbf{A}_{321})  \\
&\hspace{1cm}=(\mathbf{B}\otimes  \mathbf{I})\Big(\big(\mathbf{A}_{123}\tfrac{1}{\xi_1}+\mathbf{A}_{213}\tfrac{1}{\xi_2}\big)+ \big(\mathbf{A}_{132}\tfrac{1}{\xi_1}+\mathbf{A}_{312}\tfrac{1}{\xi_3}\big)+ \big(\mathbf{A}_{231}\tfrac{1}{\xi_2}+\mathbf{A}_{321}\tfrac{1}{\xi_3}\big)\Big).
\end{aligned}
\end{equation}
We set
\begin{equation}\label{843pm61}
\begin{aligned}
 \mathbf{A}_{123}+\mathbf{A}_{213} = ~&(\mathbf{B}\otimes  \mathbf{I})\big(\mathbf{A}_{123}\tfrac{1}{\xi_1}+\mathbf{A}_{213}\tfrac{1}{\xi_2}\big), \\[4pt]
 \mathbf{A}_{132}+\mathbf{A}_{312} = ~&(\mathbf{B}\otimes  \mathbf{I})\big(\mathbf{A}_{132}\tfrac{1}{\xi_1}+\mathbf{A}_{312}\tfrac{1}{\xi_3}\big),\\[4pt]
 \mathbf{A}_{231}+\mathbf{A}_{321} = ~&(\mathbf{B}\otimes  \mathbf{I})\big(\mathbf{A}_{231}\tfrac{1}{\xi_2}+\mathbf{A}_{321}\tfrac{1}{\xi_3}\big),
\end{aligned}
\end{equation}
to satisfy (\ref{823pm61}). Equivalently, we have
\begin{equation}\label{316am59}
\begin{aligned}
     \mathbf{A}_{213} =~& -(\mathbf{I}_{27} - \tfrac{1}{\xi_2}\mathbf{B}\otimes \mathbf{I})^{-1}(\mathbf{I}_{27} - \tfrac{1}{\xi_1}\mathbf{B}\otimes \mathbf{I})\mathbf{A}_{123},  \\[4pt]
     \mathbf{A}_{132} =~& -(\mathbf{I}_{27} - \tfrac{1}{\xi_1}\mathbf{B}\otimes  \mathbf{I})^{-1}(\mathbf{I}_{27} - \tfrac{1}{\xi_3}\mathbf{B}\otimes \mathbf{I})\mathbf{A}_{312},  \\[4pt]
     \mathbf{A}_{321} = ~&-(\mathbf{I}_{27} -\tfrac{1}{\xi_3}\mathbf{B}\otimes  \mathbf{I})^{-1}(\mathbf{I}_{27} - \tfrac{1}{\xi_2}\mathbf{B}\otimes  \mathbf{I})\mathbf{A}_{231}.
       \end{aligned}
\end{equation}
Recall that  in (\ref{316am59}) $\mathbf{I} = \mathbf{I}_3$.  Performing the matrix multiplications in (\ref{316am59}), we obtain
\begin{equation}\label{121am1117}
\begin{aligned}
 -(\mathbf{I}_{27} - \tfrac{1}{\xi_{\beta}}\mathbf{B}\otimes \mathbf{I})^{-1}(\mathbf{I}_{27} - \tfrac{1}{\xi_{\alpha}}\mathbf{B}\otimes \mathbf{I})
 ~=~&-(\mathbf{I}_{9}\otimes \mathbf{I} - \tfrac{1}{\xi_{\beta}}\mathbf{B}\otimes \mathbf{I})^{-1}(\mathbf{I}_{9}\otimes \mathbf{I} - \tfrac{1}{\xi_{\alpha}}\mathbf{B}\otimes \mathbf{I}) \\[4pt]
~=~&-\Big((\mathbf{I}_{9} - \tfrac{1}{\xi_{\beta}}\mathbf{B})\otimes \mathbf{I}\Big)^{-1}\Big((\mathbf{I}_{9} - \tfrac{1}{\xi_{\alpha}}\mathbf{B})\otimes \mathbf{I}\Big)\\[4pt]
 ~=~& -\Big((\mathbf{I}_{9} - \tfrac{1}{\xi_{\beta}}\mathbf{B})^{-1}\otimes \mathbf{I}\Big)\Big((\mathbf{I}_{9} - \tfrac{1}{\xi_{\alpha}}\mathbf{B})\otimes \mathbf{I}\Big) \\[4pt]
~=~&-\Big((\mathbf{I}_{9} - \tfrac{1}{\xi_{\beta}}\mathbf{B})^{-1}(\mathbf{I}_{9} - \tfrac{1}{\xi_{\alpha}}\mathbf{B})\Big)\otimes \mathbf{I}.
\end{aligned}
\end{equation}
Note that $\mathbf{B}$ can be made block-diagonal by reordering labels of rows and columns in the order of $11,22,33,12,21,13,31,23,32$. Using the results from (\ref{142am1019}), we then obtain
\begin{equation}\label{750pm1114}
-(\mathbf{I}_{9} - \tfrac{1}{\xi_{\beta}}\mathbf{B})^{-1}(\mathbf{I}_{9} - \tfrac{1}{\xi_{\alpha}}\mathbf{B})  = \mathbf{R}_{\beta\alpha}
\end{equation}
where $\mathbf{R}_{\beta\alpha}$ is a $9 \times 9$ matrix with the $(ij,lm)$-element given by
\begin{equation}\label{402pm1114}
\big(\mathbf{R}_{\beta\alpha}\big)_{(ij,lm)} =
\begin{cases}
S_{\beta\alpha}&~~~\textrm{if $ij = lm$ and $i\geq j$}\\[3pt]
-1&~~~\textrm{if $ij = lm$ and $i<j$}\\[3pt]
T_{\beta\alpha}&~~~\textrm{if $ij = ml$ and $i>j$}\\[3pt]
0&~~~\textrm{otherwise.}
\end{cases}
\end{equation}
Here, $S_{\beta\alpha}$ and $T_{\beta\alpha}$  are defined in (\ref{529pm625}).
Thus, (\ref{316am59}) can be  more succinctly written as
\begin{equation}\label{330am1115}
     \mathbf{A}_{\beta\alpha\gamma} =  ({\mathbf{R}}_{\beta\alpha}\otimes \mathbf{I})\mathbf{A}_{\alpha\beta\gamma}.
  \end{equation}
Therfore, if $\mathbf{A}_{\sigma}$ satisfies (\ref{330am1115}), then the Bethe ansatz solution (\ref{1118-225}) will also satisfy (\ref{408pm625}). Similarly, we can demonstrate that if
\begin{equation}\label{317am59}
         \mathbf{A}_{\alpha\gamma\beta} = (\mathbf{I} \otimes {\mathbf{R}}_{\gamma\beta})\mathbf{A}_{\alpha\beta\gamma},
 \end{equation}
then the Bethe ansatz solution (\ref{1118-225}) will satisfy (\ref{409pm625}). However, to ensure consistency between the expressions in (\ref{330am1115}) and (\ref{317am59}), it must be verified that
 \begin{equation}
  ({\mathbf{R}}_{\gamma\beta} \otimes \mathbf{I})(\mathbf{I} \otimes {\mathbf{R}}_{\gamma\alpha})({\mathbf{R}}_{\beta\alpha} \otimes \mathbf{I}) = (\mathbf{I} \otimes {\mathbf{R}}_{\beta\alpha})({\mathbf{R}}_{\gamma\alpha} \otimes \mathbf{I})(\mathbf{I} \otimes {\mathbf{R}}_{\gamma\beta})\label{134-am-529}
 \end{equation}
 and
 \begin{equation}
 ({\mathbf{R}}_{\beta\alpha} \otimes \mathbf{I})({\mathbf{R}}_{\alpha\beta} \otimes \mathbf{I}) = \mathbf{I}_{27}=( \mathbf{I}  \otimes{\mathbf{R}}_{\beta\alpha})( \mathbf{I}  \otimes{\mathbf{R}}_{\alpha\beta}).\label{135-am-529}
 \end{equation}
The proofs of (\ref{134-am-529}) and (\ref{135-am-529}) are provided in Appendix \ref{334am1115}, using the approach outlined in Section 2.4 of \cite{Lee-2021}.
\section{General $N\geq 2$}\label{1055am729}
In this section, we extend the results for $N=2,3$ to general $N\geq 2$.
\subsection{Master equations}
First, we generalize the matrices $\mathbf{B}$ in (\ref{Bmatrix}) and   (\ref{646pm623}) for $N=2,3$   to general $N$.
\begin{definition}
We define $\mathbf{B}$ as the $N^2 \times N^2$ matrix with rows and columns labelled lexicographically by $11,12,\dots, NN$,  such that the $(ij,lm)$-element is given by
\begin{equation*}
\mathbf{B}_{(ij,lm)}
= \begin{cases}
1~~~\textrm{if $ij = lm$ and $i\geq j$,}\\[5pt]
1~~~\textrm{if $ij = ml$ and $i>j$,}\\[5pt]
0~~~\textrm{otherwise},
\end{cases}
\end{equation*}
\end{definition}
Now, we define a matrix that describes the permissibility of transition for a two-particle sector consisting of the $l$-th and $(l+1)$-th particles.
\begin{definition}\label{141am731}
For $l=1,\dots, N-1$, we define $N^N \times N^N$ matrix
\begin{eqnarray}
\mathbf{B}^{(l)}&:=&\underbrace{\mathbf{I} ~~\otimes~~ \cdots ~~\otimes~~ \mathbf{I}}_{\textrm{$(l-1)$ times}} ~~\otimes~~ \mathbf{B} ~~\otimes~~ \underbrace{\mathbf{I} ~~\otimes~~ \cdots ~~\otimes~~  \mathbf{I} }_{\textrm{$(N-l-1)$ times}}; \label{510pm626}
\end{eqnarray}
where $\mathbf{I}$  is the $N \times N$ identity matrix.
\end{definition}
Then, $\mathbf{B}^{(l)}$ is an $N^N \times N^N$ matrix with rows and column labelled lexicographically  by $\pi$ and $\nu$, respectively, with $\pi, \nu = 1\cdots1, \dots, N\cdots N$. The $(\pi,\nu)$-element of $\mathbf{B}^{(l)}$ describes the permissibility of  transition
\begin{equation}\label{507pm626}
(x_1,\dots, x_{l-1},x-1, x,x_{l+2},\dots, x_N, \nu)~~~ \longrightarrow ~~~ (x_1,\dots, x_{l-1},x, x+1,x_{l+2},\dots, x_N, \pi)
\end{equation}
where $\pi_i = \nu_i$ for $i \neq  l,l+1$.

We generalize $\mathbf{B}_1$ in (\ref{Bmatrix2}) and (\ref{646pm6236}), and $\mathbf{B}_2$ in (\ref{Bmatrix3}) and (\ref{646pm62361}) to general $N$.
\begin{definition}
We define $\mathbf{B}_1$ and $\mathbf{B}_2$ as  $N^2 \times N^2$ matrices with rows and columns labelled lexicographically by $11,12,\dots, NN$,  such that their $(ij,lm)$-elements are given by
\begin{equation*}
(\mathbf{B}_1)_{(ij,lm)}
= \begin{cases}
1~~~\textrm{if $ij = ml$ and $i> j$,}\\[5pt]
0~~~\textrm{otherwise}.
\end{cases}
\end{equation*}
and
\begin{equation*}
(\mathbf{B}_2)_{(ij,lm)}
= \begin{cases}
1~~~\textrm{if $ij = lm$ and $i\geq j$,}\\[5pt]
2~~~\textrm{if $ij = lm$ and $i< j$,}\\[5pt]
0~~~\textrm{otherwise},
\end{cases}
\end{equation*}
respectively. In a same manner to Definition \ref{141am731}, we define $\mathbf{B}_1^{(l)}$ and $\mathbf{B}_2^{(l)}$.
\end{definition}
The $(\pi,\nu)$-element of $\mathbf{B}_1^{(l)}$ describes the permissibility of  transition
\begin{equation}\label{507pm626}
(x_1,\dots, x_{l-1},x, x+1,x_{l+2},\dots, x_N, \nu)~~~ \longrightarrow ~~~ (x_1,\dots, x_{l-1},x, x+1,x_{l+2},\dots, x_N, \pi)
\end{equation}
where $\pi_i = \nu_i$ for $i \neq  l,l+1$.  The $(\pi,\pi)$-element of the $N^N \times N^N$ diagonal matrix  $q\mathbf{B}_2^{(l)}$   represents  the rate of jumping out of the state
\begin{equation}\label{507pm626}
(x_1,\dots, x_{l-1},x, x+1,x_{l+2},\dots, x_N, \pi)
\end{equation}
due to the jumps of $\pi_l$ and $\pi_{l+1}$ to the left when $x-1$ is empty.

We extend the notations $\mathbf{U}(X;t)$ and $U_{(Y,\nu)}(X,\pi;t)$  used for $N=2,3$ to general $N$. If $\mathbf{U}(X;t)$ satisfies
 \begin{equation}\label{633pm1116}
 \begin{aligned}
 \frac{d}{dt}\mathbf{U}(X;t) =~& \sum_{i=1}^N p\mathbf{U}(x_1,\cdots, x_{i-1},x_i-1,x_{i+1},\cdots, x_N;t)  \\
                             +~&\sum_{i=1}^N q\mathbf{U}(x_1,\cdots, x_{i-1},x_i+1,x_{i+1},\cdots, x_N;t)\\[5pt]
                             -~&N\,\mathbf{U}(X;t)
 \end{aligned}
 \end{equation}
and the \textit{boundary conditions}
 \begin{equation}\label{717pm1116}
 \mathbf{U}(x_1,\cdots,x_{l-1},x,x,x_{l+2},\cdots, x_N;t) =\mathbf{B}^{(l)}\,\mathbf{U}(x_1,\cdots,x_{l-1},x-1,x,x_{l+2},\cdots, x_N;t)
 \end{equation}
 for all $l=1,\dots, N-1$, then $\mathbf{U}(X;t)$ satisfies the master equation corresponding to any given $X$ in the physical region. For example,  for $X=(x,x+1,\dots,x+{N-1})$, (\ref{633pm1116}) becomes
\begin{equation}\label{924pm728}
\begin{aligned}
\frac{d}{dt}\mathbf{U}(X;t) =~& p\mathbf{U}(x-1,x+1,\dots, x+N-1;t) \\
&\hspace{1cm}+~ p\mathbf{B}^{(1)}\mathbf{U}(x-1,x,x+2,\dots,x+N-1 ;t)  \\[4pt]
&\hspace{1cm}+~ p\mathbf{B}^{(2)}\mathbf{B}^{(1)}\mathbf{U}(x-1,x,x+1,x+3,\dots,x+N-1 ;t)  \\[4pt]
&\hspace{1cm}+~\cdots ~+~ p\Big(\mathbf{B}^{(N-1)}\cdots \mathbf{B}^{(1)}\Big)\mathbf{U}(x-1,x,\dots,x+N-2 ;t) \\[4pt]
&\hspace{1cm}+~q\mathbf{U}(x,x+1,\dots, x+N-2, x+N;t) \\[4pt]
&\hspace{1cm}+~ q\Big(\mathbf{B}^{(1)} + \cdots+ \mathbf{B}^{(N-1)}\Big)\mathbf{U}(X;t)- N\mathbf{U}(X;t)
\end{aligned}
\end{equation}
by the boundary conditions (\ref{717pm1116}). If we use the fact that $\mathbf{B}_2 + \mathbf{B} - \mathbf{B}_1 = 2\mathbf{I}_{N^2}$, then (\ref{924pm728}) can be rewritten as   the master equation involving  $\mathbf{B}^{(l)}, \mathbf{B}_1^{(l)}, \mathbf{B}_2^{(l)}$.
\subsection{Bethe ansatz solution}
   The Bethe ansatz solution to (\ref{633pm1116}) is given by
 \begin{equation}\label{645pm1116}
 \mathbf{U}(X;t) = \sum_{\sigma \in S_N}\mathbf{A}_{\sigma}\prod_{i=1}^N\bigg(\xi_{\sigma(i)}^{x_i} e^{\varepsilon(\xi_i)t}\bigg)
 \end{equation}
 where
 \begin{equation}\label{507am727}
 \varepsilon(\xi_i) = \frac{p}{\xi_i} +q\xi_i - 1.
 \end{equation}
 To ensure that  (\ref{645pm1116}) satisfies (\ref{717pm1116}) for each $l=1,\cdots, N-1$,   we substitute (\ref{645pm1116}) into (\ref{717pm1116}) and simplify to obtain
 \begin{equation}\label{1229am1117}
 \sum_{\sigma \in S_N}\Big(\mathbf{I}_{N^N} - \tfrac{1}{\xi_{\sigma(l)}}\mathbf{B}^{(l)}\Big) \mathbf{A}_{\sigma} = \mathbf{0}
 \end{equation}
 for each $l$.

 Let $A_N$ be the subgroup of all even permutations of ${S}_N$, and  let $T_i \in S_N$ be the simple transposition that exchanges the $i$-th and the $(i+1)$-th elements. We can then express  (\ref{1229am1117}) as
 \begin{equation}\label{1242am1117}
\mathbf{0} =  \sum_{\sigma \in A_N}\Big(\mathbf{I}_{N^N} - \tfrac{1}{\xi_{\sigma(l)}}\mathbf{B}^{(l)}\Big) \mathbf{A}_{\sigma}
 + \sum_{\sigma \in A_N}\Big(\mathbf{I}_{N^N} - \tfrac{1}{\xi_{\sigma(l+1)}}\mathbf{B}^{(l)}\Big) \mathbf{A}_{T_l\sigma}.
 \end{equation}
 To satisfy (\ref{1242am1117}), we set
 \begin{equation}\label{102pm627}
 \mathbf{0} = \Big(\mathbf{I}_{N^N} - \tfrac{1}{\xi_{\sigma(l)}}\mathbf{B}^{(l)}\Big) \mathbf{A}_{\sigma}
 + \Big(\mathbf{I}_{N^N} - \tfrac{1}{\xi_{\sigma(l+1)}}\mathbf{B}^{(l)}\Big) \mathbf{A}_{T_l\sigma}.
 \end{equation}
for each $\sigma \in A_N$. Consequently, we obtain
 \begin{equation}\label{136am1117}
 \begin{aligned}
  \mathbf{A}_{T_l\sigma} =~& -\Big(\mathbf{I}_{N^N} - \tfrac{1}{\xi_{\sigma(l+1)}}\mathbf{B}^{(1)}\Big)^{-1} ~\Big(\mathbf{I}_{N^N} - \tfrac{1}{\xi_{\sigma(l)}}\mathbf{B}^{(l)}\Big) \mathbf{A}_{\sigma} \\[5pt]
 =~& (\,\underbrace{\mathbf{I} ~\otimes~ \cdots ~\otimes~ \mathbf{I}}_{\textrm{$(l-1)$ times}}  ~\otimes~  \mathbf{R}_{\sigma(l+1)\sigma(l)} ~\otimes~ \underbrace{\mathbf{I}  ~\otimes~ \cdots ~\otimes~ \mathbf{I}}_{\textrm{$(N-l-1)$ times}}\, )\mathbf{A}_{\sigma}
 \end{aligned}
 \end{equation}
 where the second equality in (\ref{136am1117}) is derived similarly to  (\ref{121am1117}), and  $\mathbf{R}_{\sigma(l+1)\sigma(l)}$ is the $N^2 \times N^2$ matrix whose rows and columns are labelled by $11,\dots, NN$, with its $(ij,lm)$-elements given as  in (\ref{402pm1114}).
  \begin{definition}
Let $\mathbf{R}_{\beta\alpha}$ be $N^2 \times N^2$ matrix with rows and columns  labelled lexicographically by $11,\dots, NN$, where  the $(ij,lm)$-element is given by (\ref{402pm1114}). We define
 \begin{equation}\label{Rmatrix606}
 \mathbf{T}_{l,\beta\alpha}:=\,\underbrace{\mathbf{I} ~~\otimes~~ \cdots ~~\otimes~~ \mathbf{I}}_{\textrm{$(l-1)$ times}}  ~~\otimes~~  \mathbf{R}_{\beta\alpha} ~~\otimes~~ \underbrace{\mathbf{I}  ~~\otimes~~ \cdots ~~\otimes~~ \mathbf{I}}_{\textrm{$(N-l-1)$ times}}\,
 \end{equation}
where $\mathbf{I}$ is the $N\times N$ matrix.
 \end{definition}
The following  proposition follows from the arguments presented above.
\begin{proposition}
If
\begin{equation}\label{1259pm628}
 \mathbf{A}_{T_l\sigma} =  \mathbf{T}_{l,\sigma(l+1)\sigma(l)}\mathbf{A}_{\sigma}
 \end{equation}
for each $\sigma  \in S_N$, then the boundary condition (\ref{717pm1116}) is satisfied.
\end{proposition}
Since $\{T_1,\cdots, T_{N-1}\}$ generates $S_N$, for any $\sigma \in S_N$, there exists an expression
\begin{equation}\label{112pm628}
\sigma = T_{l_k}\cdots T_{l_2}T_{l_1}
\end{equation}
for some $l_1,\dots, l_k \in \{1,\dots, N-1\}$.  The following lemma provided a condition for (\ref{1259pm628}) to be satisfied, which is evident  from the construction.
\begin{lemma}\label{347pm628}
Let $\sigma$ be given by (\ref{112pm628}) and let $\sigma^{(i)} =  T_{l_i}\cdots T_{l_1}$, $i=1,\dots, k$ with $\sigma^{(0)}$ denoting the identity permutation. If
\begin{equation}\label{150am727}
\mathbf{A}_{\sigma} = \mathbf{A}_{ T_{l_k}\cdots T_{l_2}T_{l_1}} := \mathbf{T}_{l_k,\sigma^{(k-1)}(l_k+1)\sigma^{(k-1)}(l_k)}~ \times ~\cdots~ \times~ \mathbf{T}_{l_1,\sigma^{(0)}(l_1+1)\sigma^{(0)}(l_1)},
\end{equation}
then (\ref{1259pm628}) is satisfied.
\end{lemma}
For $\mathbf{A}_{\sigma}$ in (\ref{150am727}) to be well-defined,  the following relations  must also be satisfied:
\begin{equation}\label{416pm628}
\begin{aligned}
  &\mathbf{T}_{i,\alpha\beta}\mathbf{T}_{i,\beta\alpha} = \mathbf{I};\\
   &\mathbf{T}_{i,\beta\gamma}\mathbf{T}_{i+1,\alpha\gamma}\mathbf{T}_{i,\alpha\beta} = \mathbf{T}_{i+1,\alpha\beta}\mathbf{T}_{i,\alpha\gamma}\mathbf{T}_{i+1,\beta\gamma}.
\end{aligned}
\end{equation}
In fact, Lemma \ref{1242pm628} and Lemma \ref{417pm628} demonstrate that these relations are satisfied.
\subsection{Transition probabilities}\label{109-2023-648pm}
In the previous section, we claimed that $\mathbf{U}(X;t)$ in (\ref{645pm1116}) with $\mathbf{A}_{\sigma}$ defined by (\ref{150am727}) satisfies the master equation for each $X$ in the corresponding physical region. This holds true when performing contour integrals over all $\xi_i$ variables with appropriate contours. This approach is well known for various integrable particle models \cite{Schutz-1997,Tracy-Widom-2008}. The  transition probability $P_{(Y,\nu)}(X,\pi;t)$ of our model can also be expressed as a contour integral, analogous to the formula for the multi-species ASEP in \cite{Lee-2020}.
\begin{theorem}
Let $C_{R_i}$ be a counterclockwise circle centered at the origin with radius $R_i$ and suppose that $1<R_1< \cdots <R_N$. Let $\mathbf{A}_{\sigma}$ be given by (\ref{150am727}) and $\varepsilon(\xi_i)$ by (\ref{507am727}). For any arbitrary initial state $(Y,\nu)$, the transition probability $P_{(Y,\nu)}(X,\pi;t)$ is given by:
\begin{equation}\label{520pm727}
P_{(Y,\nu)}(X,\pi;t) = \Big(\frac{1}{2\pi i}\Big)^N\int_{C_{R_N}}\cdots \int_{C_{R_1}} \sum_{\sigma \in S_N}\big(\mathbf{A}_{\sigma}\big)_{(\pi,\nu)}\Big(\prod_{i=1}^N\xi_{\sigma(i)}^{x_i - y_{\sigma(i)} - 1}e^{\varepsilon(\xi_i)t}\Big)d\xi_1 \dots d\xi_N.
\end{equation}
\end{theorem}
\begin{proof}
It is sufficient to verify the initial condition. Since the denominators of the nontrivial terms in the matrix $\mathbf{R}_{\beta\alpha}$ that is contained in the expression of $\mathbf{A}_{\sigma}$ match those of $S_{\beta\alpha}^{\dagger}$ with $\mu=1$ and $\lambda = 0$ in Section 4 of \cite{Lee-2012}, the proofs of Corollary 2 and Theorem 1 in \cite{Lee-2012} essentially imply that (\ref{520pm727}) satisfies the initial condition.
\end{proof}
\section{Discussion}
We have defined a multi-species stochastic particle model (distinct from the multi-species ASEP) on $\mathbb{Z}$, where particles move either  to the right or the left: a particle jumping to the right moves to the nearest site occupied by a lower-numbered particle (or an empty site) and pushes the lower-numbered particle  forward (if any);  a particle jumping  to the left follows the rules of the conventional multi-species TASEP. We have shown that  this model is integrable and  derived the formula for the transition probability for an arbitrary initial state. A key motivation of defining this model such a way that the model can be integrable is the symmetry between the one-sided model with long-range jumps and the multi-species TASEP in the opposite direction, in the sense that their boundary conditions are compatible when using the Bethe ansatz.

Given the form of the matrices in (\ref{520pm727}), it is anticipated that $P_{(Y,\nu)}(X,\pi;t)$ with $\pi = \nu$ can be expressed as a determinant (See \cite{Lee-2020} for related results in the multi-species TASEP context). Generally, the integrand in (\ref{520pm727}) can be a sum of numerous products due to  $\mathbf{A}_{\sigma}$ being a product of matrices as described in (\ref{150am727}).  However, for certain special initial permutations $\nu$ of species, such as $12\cdots N$, $12\cdots2$, $1\cdots 12$ and so on, it seems that $\big(\mathbf{A}_{\sigma}\big)_{(\pi,\nu)}$ can be expressed as a single product  (See \cite{Lee-Raimbekov-2022} for related results for the multi-species ASEP). Also, it would be interesting to investigate whether the methods used in this paper could serve as a machinery to develop other new integrable multi-species particle models as well as multi-species versions of other known single-species models such as \cite{Povolotsky-2}. Furthermore, these methods might help find the transition probabilities for some known multi-species particle models mentioned in \cite{Kuan-2018,Kuniba-Okado-Watanabe,Takeyama}.

Finally, we have not yet identified an integrable multi-species model where particles make long-range jumps in both directions, that is, a multi-species version of the model in \cite{Ali2}. Developing such a multi-species version and finding its transition probability remain an area for future research.

\begin{appendix}
\section{Proof of the consistency conditions}\label{334am1115}
We use the same approach as in Section 2.4 of \cite{Lee-2021}.
\begin{lemma}\label{1242pm628}
If $\mathbf{R}_{\beta\alpha}$ is the $N^2 \times N^2$ matrix defined by (\ref{402pm1114}), then  $\mathbf{R}_{\beta\alpha}\mathbf{R}_{\alpha\beta}$ is the identity matrix.
\end{lemma}
\begin{proof}
Let $ij$ and $kl$ be labels of rows and columns with $ij, kl = 11,\dots, NN$. For $i =j$,
 \begin{equation*}
 \begin{aligned}
 \big(\mathbf{R}_{\beta\alpha}\mathbf{R}_{\alpha\beta}\big)_{(ii,kl)} =~& \sum_{(mn)}\big(\mathbf{R}_{\beta\alpha}\big)_{(ii,mn)}\big(\mathbf{R}_{\alpha\beta}\big)_{(mn,kl)} \\
 =~&\big(\mathbf{R}_{\beta\alpha}\big)_{(ii,ii)}\big(\mathbf{R}_{\alpha\beta}\big)_{(ii,kl)} = \begin{cases}
 S_{\beta\alpha}S_{\alpha\beta}=1~&\textrm{if $kl = ii$},\\[4pt]
 S_{\beta\alpha}\cdot 0=0~&\textrm{if $kl \neq ii$}.
 \end{cases}
 \end{aligned}
 \end{equation*}
For $i>j$,
\begin{equation*}
 \begin{aligned}
 \big(\mathbf{R}_{\beta\alpha}\mathbf{R}_{\alpha\beta}\big)_{(ij,kl)} =~& \sum_{(mn)}\big(\mathbf{R}_{\beta\alpha}\big)_{(ij,mn)}\big(\mathbf{R}_{\alpha\beta}\big)_{(mn,kl)} \\[4pt]
 =~&\big(\mathbf{R}_{\beta\alpha}\big)_{(ij,ij)}\big(\mathbf{R}_{\alpha\beta}\big)_{(ij,kl)} + \big(\mathbf{R}_{\beta\alpha}\big)_{(ij,ji)}\big(\mathbf{R}_{\alpha\beta}\big)_{(ji,kl)}\\[4pt]
 =~&\begin{cases}
  S_{\beta\alpha}\cdot S_{\alpha\beta} + T_{\beta\alpha}\cdot 0 =1 ~&\textrm{if $kl = ij$},\\[4pt]
 S_{\beta\alpha}\cdot T_{\alpha\beta}+T_{\beta\alpha}\cdot(-1) = 0~&\textrm{if $kl = ji$} \\[4pt]
 S_{\beta\alpha}\cdot0 + T_{\beta\alpha}\cdot 0 = 0~& \textrm{if $kl \neq ij,ji$}.
 \end{cases}
 \end{aligned}
 \end{equation*}
 For $i<j$,
 \begin{equation*}
 \begin{aligned}
 \big(\mathbf{R}_{\beta\alpha}\mathbf{R}_{\alpha\beta}\big)_{(ij,kl)} =~& \sum_{(mn)}\big(\mathbf{R}_{\beta\alpha}\big)_{(ij,mn)}\big(\mathbf{R}_{\alpha\beta}\big)_{(mn,kl)} \\[4pt]
 =~&\big(\mathbf{R}_{\beta\alpha}\big)_{(ij,ij)}\big(\mathbf{R}_{\alpha\beta}\big)_{(ij,kl)}+\big(\mathbf{R}_{\beta\alpha}\big)_{(ij,ji)}\big(\mathbf{R}_{\alpha\beta}\big)_{(ji,kl)}\\[4pt]
 =~&\begin{cases}
  (-1)\cdot(-1)+  0\cdot T_{\beta\alpha} =1 ~&\textrm{if $kl = ij$},\\[4pt]
  (-1)\cdot 0+0\cdot S_{\alpha\beta} = 0~&\textrm{if $kl = ji$} \\[4pt]
 (-1)\cdot 0 +0\cdot  0= 0~&\textrm{if $kl \neq ij,ji$}.
 \end{cases}
 \end{aligned}
 \end{equation*}
\end{proof}
\begin{lemma}\label{417pm628} (Yang-Baxter equation)
Let $N \geq 3$. If $\mathbf{R}_{\beta\alpha}$ is the $N^2 \times N^2$ matrix defined by (\ref{402pm1114}) and $\mathbf{I}$ be the $N \times N$ matrix,
\begin{equation}\label{134-am-52976}
({\mathbf{R}}_{\gamma\beta} \otimes \mathbf{I})(\mathbf{I} \otimes {\mathbf{R}}_{\gamma\alpha})({\mathbf{R}}_{\beta\alpha} \otimes \mathbf{I}) = (\mathbf{I} \otimes {\mathbf{R}}_{\beta\alpha})({\mathbf{R}}_{\gamma\alpha} \otimes \mathbf{I})(\mathbf{I} \otimes {\mathbf{R}}_{\gamma\beta}).
\end{equation}
\end{lemma}
\begin{proof}
Both $(\mathbf{I} \otimes {\mathbf{R}}_{\beta\alpha})$ and $({\mathbf{R}}_{\beta\alpha} \otimes \mathbf{I})$ are $N^3 \times N^3$ matrices with rows and columns labelled with $111,\cdots, NNN$, and each label is a permutation of a certain multi-set $[i,j,k]$. If $ijk$ and $lmn$ are permutations from two different multi-sets, then
\begin{equation*}
(\mathbf{I} \otimes {\mathbf{R}}_{\beta\alpha})_{(ijk,lmn)} = 0 = ({\mathbf{R}}_{\beta\alpha} \otimes \mathbf{I})_{(ijk,lmn)}.
\end{equation*}
 Hence, if we reorder the columns and rows so that all permutations  of a given multi-set are grouped, then both $(\mathbf{I} \otimes {\mathbf{R}}_{\beta\alpha})$ and $({\mathbf{R}}_{\beta\alpha} \otimes \mathbf{I})$ become block-diagonal. Let
\begin{equation*}
({\mathbf{R}}_{\beta\alpha} \otimes \mathbf{I})_{[i,j,k]}
\end{equation*}
be the block on the diagonal of the matrix ${\mathbf{R}}_{\beta\alpha} \otimes \mathbf{I}$  whose  rows are columns are labelled with permutations of $[i,j,k]$. We similarly define  $(\mathbf{I} \otimes {\mathbf{R}}_{\beta\alpha})_{[i,j,k]}$. It suffices to show  (\ref{134-am-52976}) \textit{blockwise}, that is,
\begin{equation}\label{357pm1115}
({\mathbf{R}}_{\gamma\beta} \otimes \mathbf{I})_{[i,j,k]}(\mathbf{I} \otimes {\mathbf{R}}_{\gamma\alpha})_{[i,j,k]}({\mathbf{R}}_{\beta\alpha} \otimes \mathbf{I})_{[i,j,k]} = (\mathbf{I} \otimes {\mathbf{R}}_{\beta\alpha})_{[i,j,k]}({\mathbf{R}}_{\gamma\alpha} \otimes \mathbf{I})_{[i,j,k]}(\mathbf{I} \otimes {\mathbf{R}}_{\gamma\beta})_{[i,j,k]}
\end{equation}
for each multi-set $[i,j,k]$.  Each multi-set $[i,j,k]$ of integers $1,\dots, N$ is one of the forms $[i,i,i]$, $[i,i,j]$ with $i< j$, $[i,i,j]$ with $i> j$ and $[i,j,k]$ with $i<j<k$.  First, for the multi-set $[i,i,i]$, (\ref{357pm1115}) is simply
\begin{equation*}
S_{\gamma\beta}S_{\gamma\alpha}S_{\beta\alpha}  = S_{\beta\alpha}S_{\gamma\alpha}S_{\gamma\beta}
\end{equation*}
which is obviously true because
\begin{equation*}\label{532pm1115}
({\mathbf{R}}_{\beta\alpha} \otimes \mathbf{I})_{[i,i,i]} = S_{\beta\alpha} = (\mathbf{I} \otimes {\mathbf{R}}_{\beta\alpha})_{[i,i,i]}.
\end{equation*}
For $[i,i,j]$ with $i<j$, we can obtain
\begin{equation*}\label{533pm1115}
({\mathbf{R}}_{\beta\alpha} \otimes \mathbf{I})_{[i,i,j]}  =
\kbordermatrix{
 & iij & iji & jii \\
iij & S_{\beta\alpha}& 0 & 0  \\
iji &0 &-1 & 0\\
jii & 0& T_{\beta\alpha} & S_{\beta\alpha}
}
\end{equation*}
and
\begin{equation*}\label{533pm11159}
(\mathbf{I} \otimes {\mathbf{R}}_{\beta\alpha})_{[i,i,j]}=\kbordermatrix{
 & iij & iji & jii \\
iij & -1& 0 & 0  \\
iji &T_{\beta\alpha} &S_{\beta\alpha} & 0\\
jii & 0& 0 & S_{\beta\alpha}.
}
\end{equation*}
Using these matrices, we can verify (\ref{357pm1115}) by direct matrix computation. Similarly, (\ref{357pm1115}) is verified for $[i,i,j]$ with $i>j$. Finally, for $[i,j,k]$ with $i<j<k$, we obtain
\begin{equation*}\label{535pm1115}
({\mathbf{R}}_{\beta\alpha} \otimes \mathbf{I})_{[i,j,k]}~ =~ \kbordermatrix{
     & ijk & ikj & jik & jki & kij & kji \\
    ijk & -1 & 0 & 0 & 0 & 0 & 0 \\
    ikj & 0 & -1 & 0 & 0 & 0 & 0 \\
   jik &  T_{\beta\alpha} & 0 & S_{\beta\alpha} & 0 &0  & 0 \\
   jki & 0 & 0 & 0 & -1 & 0 &  0\\
   kij & 0 &  T_{\beta\alpha} & 0 & 0 &S_{\beta\alpha}  & 0 \\
   kji & 0 & 0 & 0 & T_{\beta\alpha} & 0 & S_{\beta\alpha} \\
   }
\end{equation*}
and
\begin{equation*}\label{536pm1115}
(\mathbf{I} \otimes {\mathbf{R}}_{\beta\alpha})_{[i,j,k]}~ =~ \kbordermatrix{
     & ijk & ikj & jik & jki & kij & kji \\
    ijk & -1 & 0 & 0 & 0 & 0 & 0 \\
    ikj & T_{\beta\alpha} & S_{\beta\alpha} & 0 & 0 & 0 & 0 \\
   jik & 0  & 0 & -1 & 0 &0  & 0 \\
   jki & 0 & 0 & T_{\beta\alpha} & S_{\beta\alpha} & 0 &  0\\
   kij & 0 &  0 & 0 & 0 &-1  & 0 \\
   kji & 0 & 0 & 0 & 0 &  T_{\beta\alpha} & S_{\beta\alpha} \\
   }.
\end{equation*}
Again, we can verify (\ref{357pm1115}) by  directly performing matrix multiplication.
\end{proof}

\end{appendix}

\textbf{Acknowledgement} This research is funded by Nazarbayev University under the faculty-development competitive research grants program
for 2024 -- 2026 (grant number 021220FD4251). The author is grateful to Axel R. Saenz for valuable discussions, and to Temirlan Raimbekov and Nazgul Tileukabyl for their assistance in writing the manuscript.
\\ \\
\textbf{Conflict of Interest} The author declares that there is no conflict of interest.

\end{document}